\documentclass[12pt]{amsart}
\usepackage[utf8]{inputenc}
\usepackage{graphicx}
\usepackage{xcolor}
\usepackage{times}
\usepackage[hidelinks]{hyperref}
\usepackage{enumerate,latexsym}
\usepackage{latexsym}
\usepackage{amssymb}
\usepackage{float}

\usepackage{comment}
\usepackage{amsfonts,amssymb}
\usepackage{dsfont}

\makeatletter
\def\namedlabel#1#2{\begingroup
 #2%
 \def\@currentlabel{#2}%
 \phantomsection\label{#1}\endgroup
}
\makeatother

\newtheorem{theorem}{Theorem}[section]
\newtheorem*{acknowledgement*}{Acknowledgements}

\newtheorem{conjecture}[theorem]{Conjecture}
\newtheorem{corollary}[theorem]{Corollary}

\newtheorem{question}[theorem]{Question}
\newtheorem{proposition}[theorem]{Proposition}
\theoremstyle{definition}
\newtheorem{definition}[theorem]{Definition}

\theoremstyle{definition}

\newcommand{\ben}{\begin{enumerate}}
\newcommand{\een}{\end{enumerate}}

\newcommand{\vol}{\mbox{\rm Vol}}

\newcommand{\ed}{\end{document}}

\definecolor{rrr}{rgb}{.9,0,.1}

\definecolor{rr}{rgb}{.8,0,.3}
\definecolor{pp}{rgb}{.5,0,.7}

\usepackage{pdfsync}

\title{Hyperbolic Knotoids}
\author[C. Adams et al]{Colin Adams}
\address{Department of Mathematics, Williams College, Williamstown, MA 01267}
\email{cadams@williams.edu}
\author[]{Alexandra Bonat}\address{Department of Mathematics, Williams College, Williamstown, MA 01267} \email{ajb10@williams.edu}
\author[]{Maya Chande}\address{Department of Mathematics, Princeton University, Princeton, NJ 08544}
\email{mchande@princeton.edu}
\author[]{Joye Chen} \address{Department of Mathematics, Princeton University, Princeton, NJ 08544} \email{joyec@princeton.edu}
\author[]{Maxwell Jiang}
\address{Department of Mathematics, M.I.T., Cambridge, MA 02139}
\email{mdjiangt@mit.edu}
\author[]{Zachary Romrell}
\address{Department of Mathematics, Williams College, Williamstown, MA 01267} \email{zr3@williams.edu}
\author[]{Daniel Santiago}
\address{Department of Mathematics, M.I.T., Cambridge, MA 02139}
\email{dsantiag@mit.edu}
\author[]{Benjamin Shapiro}\address{Department of Mathematics, Williams College, Williamstown, MA 01267} \email{bis1@williams.edu}
\author[]{Dora Woodruff}
\address{Department of Mathematics, Harvard University, Cambridge, MA 02138}
\email{dorawoodruff@college.harvard.edu}


\begin{document}  

\begin{abstract} In  2010, Turaev introduced knotoids as a variation on knots that replaces the embedding of a circle with the embedding of a closed interval with two endpoints. A variety of knot invariants have been extended to knotoids. Here we provide definitions of hyperbolicity for both spherical and planar knotoids. We prove that the product of hyperbolic spherical knotoids is hyperbolic and the volumes add. We also determine the least volume of a rational spherical knotoid and provide various classes of hyperbolic knotoids. We also include tables of hyperbolic volumes for both spherical and planar knotoids.  
\end{abstract}

\maketitle




\section{Introduction}

\begin{figure}[htbp]
    \centering
    \includegraphics[scale=0.3]{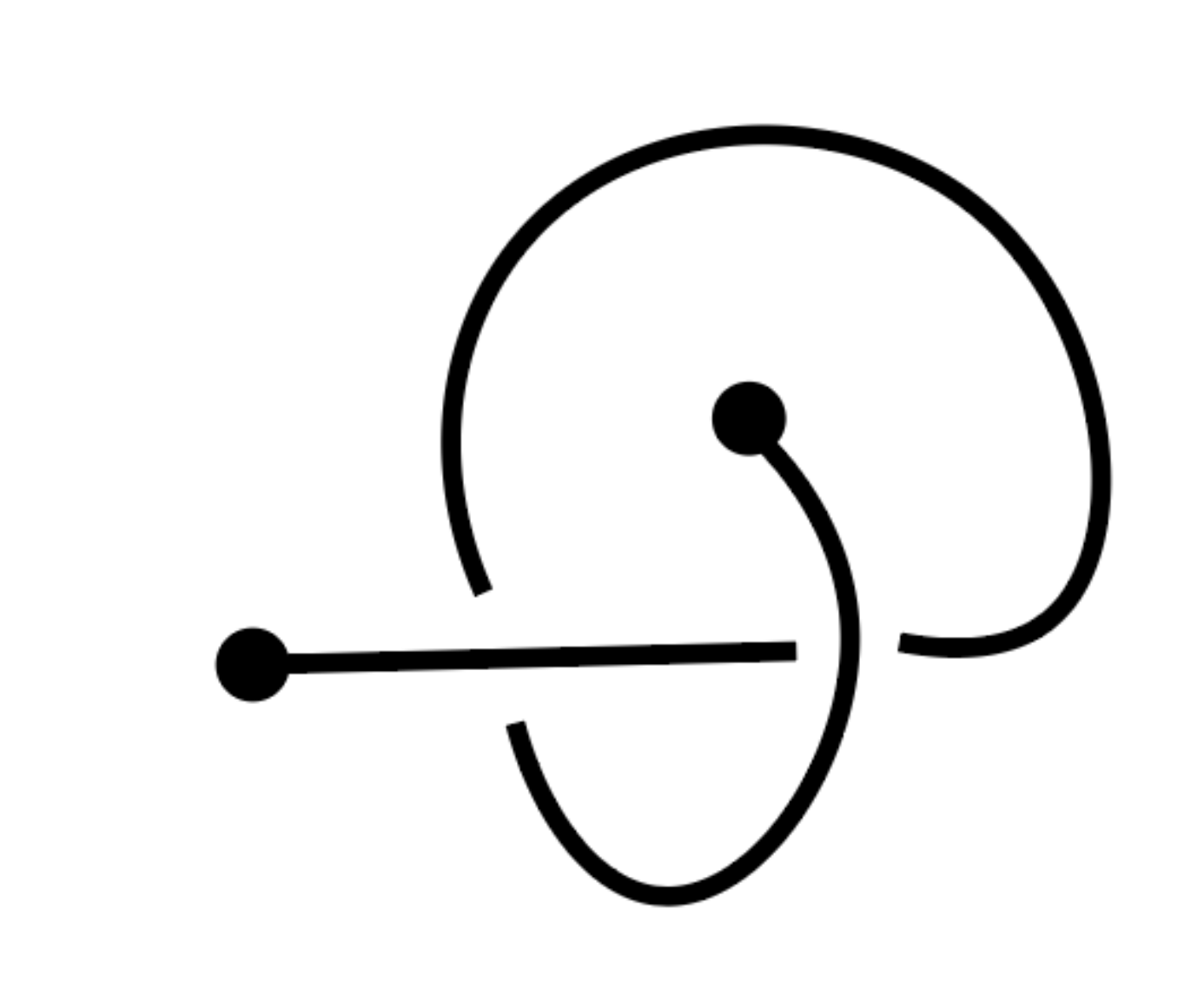}
    \caption{An example of a knotoid diagram. Throughout the paper, we refer to this knotoid as $2_1$, following the notation of \cite{tabulation}.}
    \label{knotoid_projection}
\end{figure}

In \cite{turaev}, Turaev introduced a generalization of knots called \textit{knotoids}. A \textit{knotoid diagram} is a closed interval $I = [0, 1]$ immersed in a surface $\Sigma$ with under/overcrossing data at the double points, as appears in Figure \ref{knotoid_projection}. We sometimes refer to $\Sigma$ as the \textit{projection surface} for the diagram. 

\begin{figure}[htbp]
\centering
\includegraphics[scale=0.4]{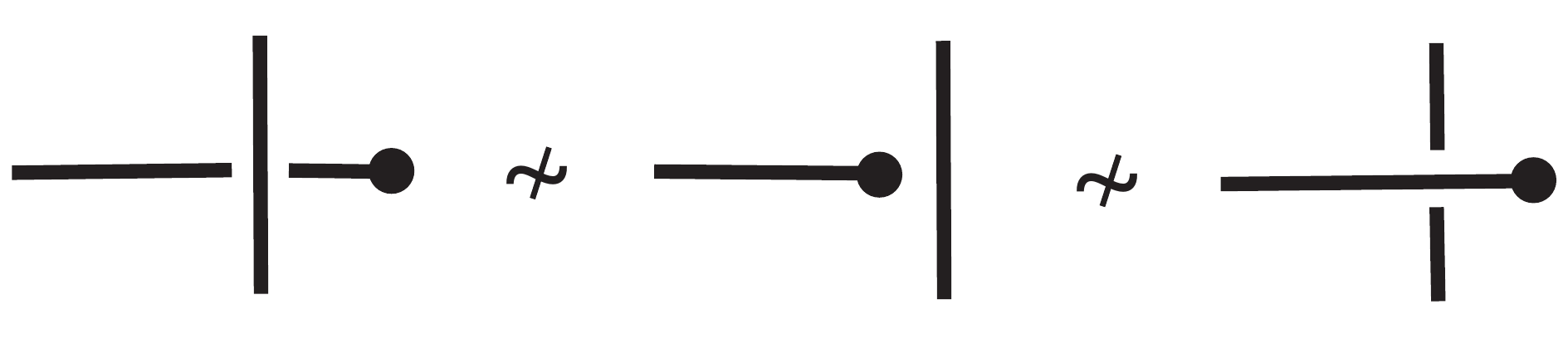}
\caption{We are forbidden from passing an endpoint over or under a strand. }
\label{forbiddenmove}
\end{figure}

 We say two knotoid diagrams are equivalent if and only if we can obtain one from the other using a sequence of Reidemeister moves away from the endpoints. We do not allow the fourth \textit{forbidden} Reidemeister move as in Figure \ref{forbiddenmove} (otherwise, we could just shrink one endpoint of the knotoid along the strand until we obtain a segment with no crossings, and the theory would be trivial). Then a knotoid is the equivalence class of a knotoid diagram up to these moves. 
 
\begin{figure}[htbp]
\centering
\includegraphics[scale=0.6]{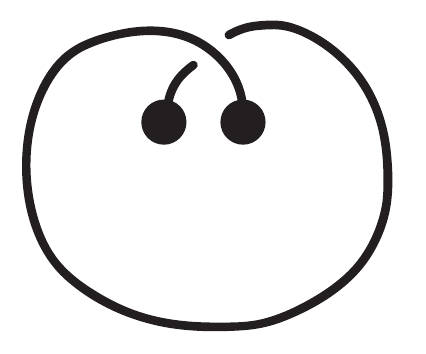}
\caption{This knotoid, denoted $1_1$, is nontrivial when viewed in the plane but trivial when viewed in $S^2$, since we can pass the strand over the back of the sphere. It is also an example of a knot-type knotoid.}
\label{1_1}
\end{figure}

 In this paper, we will usually consider knotoids in the sphere or in the plane, and we denote these sets of knotoids by $\mathbb{K}(S^2)$ and $\mathbb{K}(\mathbb{R}^2)$ respectively. The distinction between these projection surfaces is significant and gives rise to different theories: there are pairs of knotoid diagrams which are equivalent in the sphere but inequivalent in the plane. See Figure \ref{1_1} for an example. More generally, it is possible to consider $\mathbb{K}(\Sigma)$ for any projection surface $\Sigma$. Note that the theory of knotoids in $\mathbb{R}^2$ is equivalent to the theory of knotoids in $D^2$, so we interchange between these projection surfaces as needed.

In Turaev's original paper \cite{turaev} and many subsequent papers, authors have generalized various classical knot invariants to knotoids. In this paper, the main goal is to associate to a spherical or planar knotoid a well-defined hyperbolic volume and study this new invariant. Hyperbolic volume has been considered as an invariant in \cite{BBHL} using double branched covers, however, we use different constructions to define hyperbolic volume as described below.

In Section \ref{spherical_knotoid_section}, we describe two ways to map spherical knotoids to knots in the thickened torus, one of which was first defined by Turaev in \cite{turaev}. Then we may extend the notion of hyperbolicity to $\mathbb{K}(S^2)$ by considering the hyperbolicity of the image of the knotoid under these maps. We show that both of our constructions generate the same volume. 

One particularly nice feature is that the product of two hyperbolic spherical knotoids is also hyperbolic, and moreover, the hyperbolic volume is additive under the operation of taking products of spherical knotoids. We also prove other fundamental properties of this invariant and exhibit a large class of hyperbolic spherical knotoids.


In Section \ref{vol_bounds_spherical}, we prove volume bounds for hyperbolic families of spherical knotoids. One of our main results is to determine the unique rational knotoid of minimum volume. 

In Section \ref{planar_knotoid_section}, we describe two maps from planar knotoids to knots in handlebodies. Then we define a planar knotoid as hyperbolic if its image is tg-hyperbolic, meaning the complement of the image knot in the handlebody admits a hyperbolic metric such that the higher genus boundary components are totally geodesic in the metric. Then we can define the hyperbolic volume of planar knotoids and study some fundamental properties of this new invariant. We exhibit an infinite family of hyperbolic planar knotoids using results from \cite{CompositionPaper}. 
Finally, in Section \ref{volume_tables}, we provide tables of volume computations for spherical and planar knotoids. \\

For the remainder of this section, we recall some definitions that will be useful to us throughout the paper. We can consider knotoids with or without orientation (when not specified, we assume knotoids are unoriented). For an oriented knotoid in the sphere or plane, its {\it tail} and {\it head} are the endpoints corresponding to the endpoints $\{0\}$ and $\{1\}$ from the interval $I$ respectively. Given a knotoid diagram $\mathcal{D}$, a \emph{shortcut} is an arc connecting the endpoints of the diagram that does not pass through any crossings and intersects $\mathcal{D}$ transversely. Then for a knotoid $k$, its \emph{height $h(k)$} is the minimum number of intersection points between any shortcut and $\mathcal{D}$, not including the endpoints, minimized over all diagrams $\mathcal{D}$. A knotoid of height zero is called \emph{knot-type}. For instance, the planar knotoid $1_1$ is knot-type (see Figure \ref{1_1}). A knotoid that is not knot-type is said to be \emph{proper}.

\begin{figure}[htbp]
\centering
\includegraphics[scale=0.5]{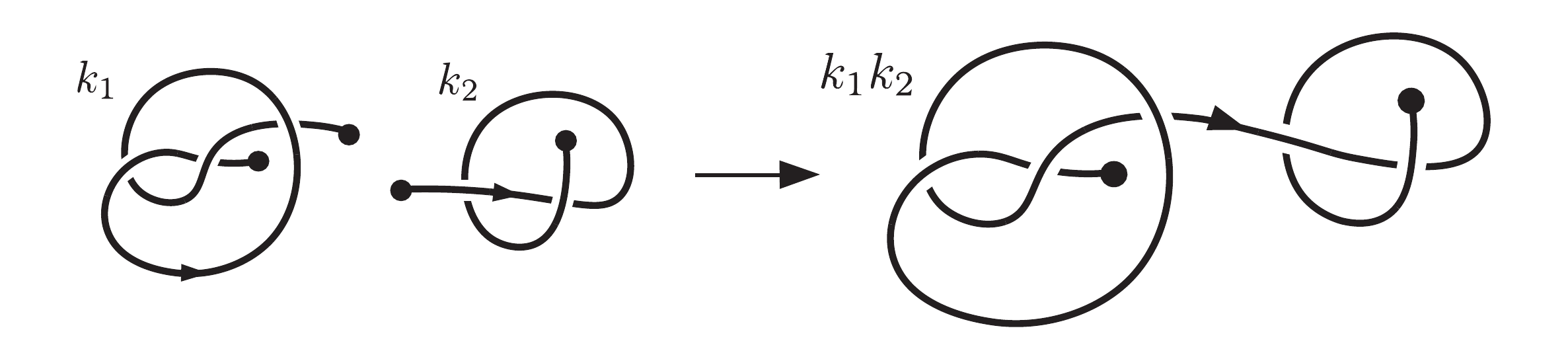}
\caption{The product of two oriented knotoids, which can be viewed on the sphere or the plane.}
\label{product}
\end{figure}

Given two oriented spherical knotoids $k_1, k_2 \in \mathbb{K}(S^2)$, we can define their \textit{product} $k_1k_2$ to be the spherical knotoid obtained by concatenating the head of $k_1$ with the tail of $k_2$, as shown in Figure \ref{product}. 

Similarly, we may also define the \emph{product} $k_1k_2$ of two oriented planar knotoids $k_1, k_2 \in \mathbb{K}(\mathbb{R}^2)$ by concatenating the head of $k_1$ with the tail of $k_2$; however, this is only possible when the tail of $k_2$ is in the exterior region (meaning the unbounded complementary region of $\Sigma \setminus \mathcal{D}$) of its projection surface $\Sigma = \mathbb{R}^2$. If the head of $k_1$ is not in the exterior region, we shrink $k_2$ to fit in the region of $k_1$ that contains the head of $k_1$ and then attach the tail of $k_2$ to the head of $k_1$.

\begin{figure}[htbp]
\centering
\includegraphics[scale=0.6]{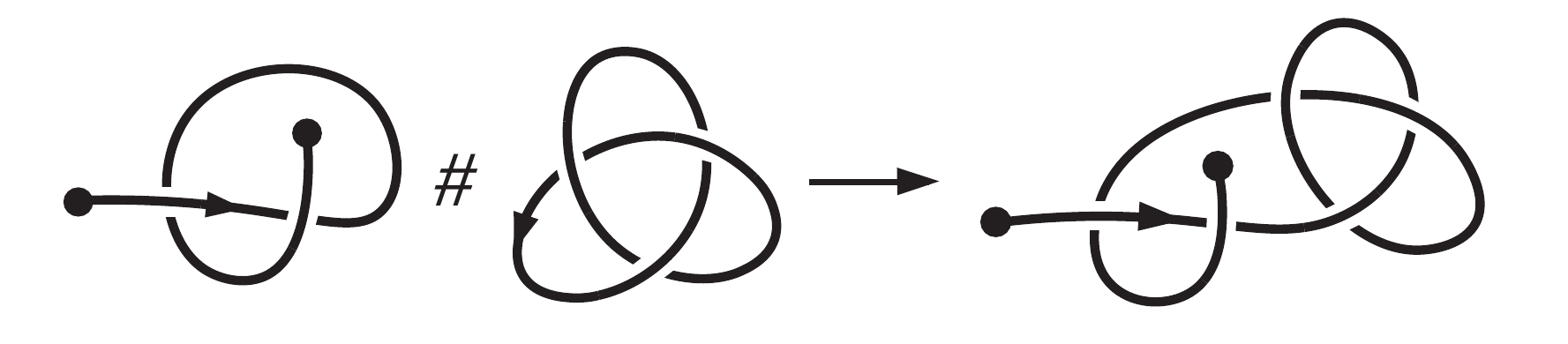}
\caption{Composition of an oriented knotoid with a knot. This can be viewed on the sphere or in the plane.}
\label{connected_sum}
\end{figure}

We may also define the \textit{composition} of an oriented knotoid $k \in \mathbb{K}(\Sigma)$ for $\Sigma = S^2$ or $\mathbb{R}^2$ with an oriented knot $K \in \mathcal{K}(\Sigma)$, by removing an arc each from $k$ and $K$, then identifying the endpoints in an orientation-preserving way. The result is an oriented knotoid in $\Sigma$ which we denote by $k \# K$.  See Figure \ref{connected_sum} for an example.

A knotoid is \emph{prime} if it is not expressible as the product of two knotoids, and it is \emph{knot-free} if it is not expressible as the composition of a knotoid with a knot. We will see that being knot-free is a necessary condition for hyperbolicity. 

The notion of \textit{rail diagrams}, first defined and studied in \cite{rails}, will also be useful to us. For any orientable, compact surface, consider the thickened surface $\Sigma \times I$ together with two distinct line segments $\ell_1 = \{x\} \times I, \ell_2 = \{y\} \times I$. Then a \textit{rail diagram} (on $\Sigma$) is a proper embedding of an arc into $(\Sigma \times I) \setminus \mathring{N}(\ell_1 \cup \ell_2)$ such that one endpoint is embedded in $\partial N(\ell_1)$ and one is embedded in $\partial N(\ell_2)$. We consider rail diagrams up to proper isotopies of the embedded arc in $(\Sigma \times I) \setminus \mathring{N}(\ell_1 \cup \ell_2)$. Note that the arc is forbidden from passing through a rail, and in this sense, the rail diagram captures topologically the forbidden move for knotoid diagrams on $\Sigma$. Indeed, the authors of \cite{rails} prove that in the planar case $\Sigma = D^2$, the theory of rail diagrams is equivalent to the theory of planar knotoids. In fact, their proof goes through for any surface $\Sigma$. Hence, we sometimes refer to both the embedded arc in the rail diagram and the corresponding knotoid by $k$. 

\begin{center}
\begin{figure}[htbp]
\includegraphics[scale=0.7]{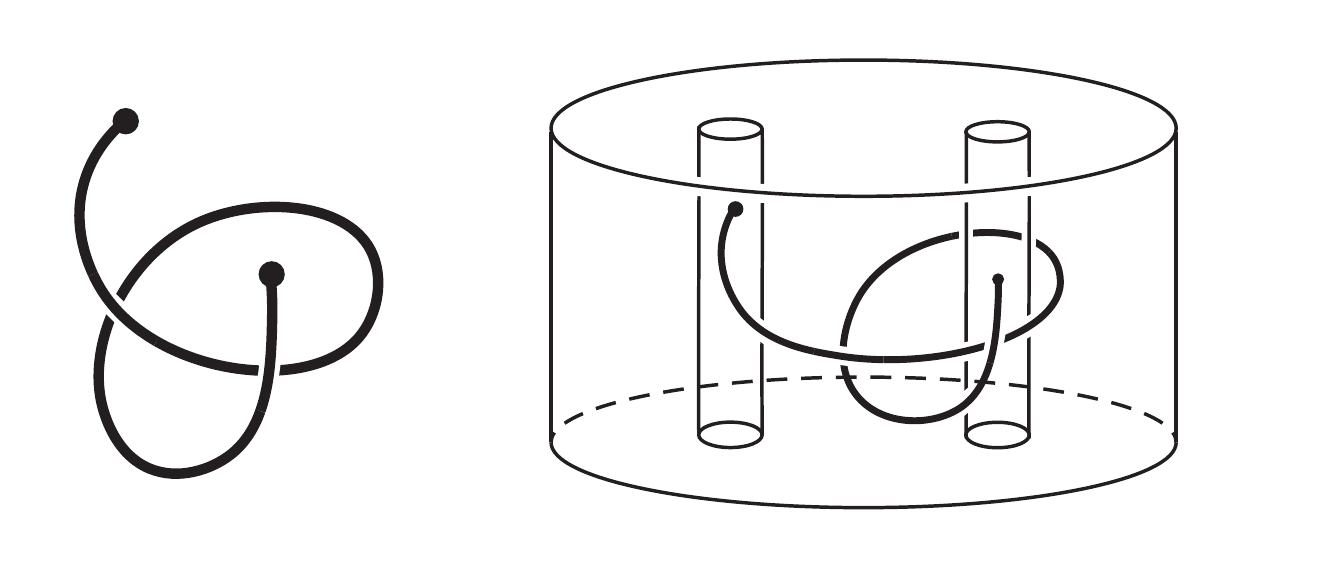}
\caption{The planar knotoid $2_1$ and its rail diagram.}
\label{rail_diagram}
\end{figure}
\end{center}

Given any knotoid in any projection surface, there are two naturally associated knots: the \emph{underclosure} and the \emph{overclosure}. The underclosure is defined by adding a shortcut to the knotoid, and having the shortcut pass under any strand it intersects; the overclosure is defined similarly. The overclosure and underclosure are equivalent for height zero knotoids. 

Another type of knotoid closure is the \textit{virtual closure}. Kauffman defined virtual knot theory in \cite{Virtual_knots}; this is a generalization of classical knot theory introducing virtual crossings in addition to the usual under and overcrossings. Virtual knots are considered equivalent up to classical Reidmeister moves, virtual Reidemeister moves, and mixed Reidemeister moves. The virtual closure of a knotoid is defined by adding a shortcut to the knotoid and making any crossings between the shortcut and knotoid into virtual crossings, yielding a virtual knot.

\begin{center}
\begin{figure} [htbp]
\includegraphics[scale=0.4]{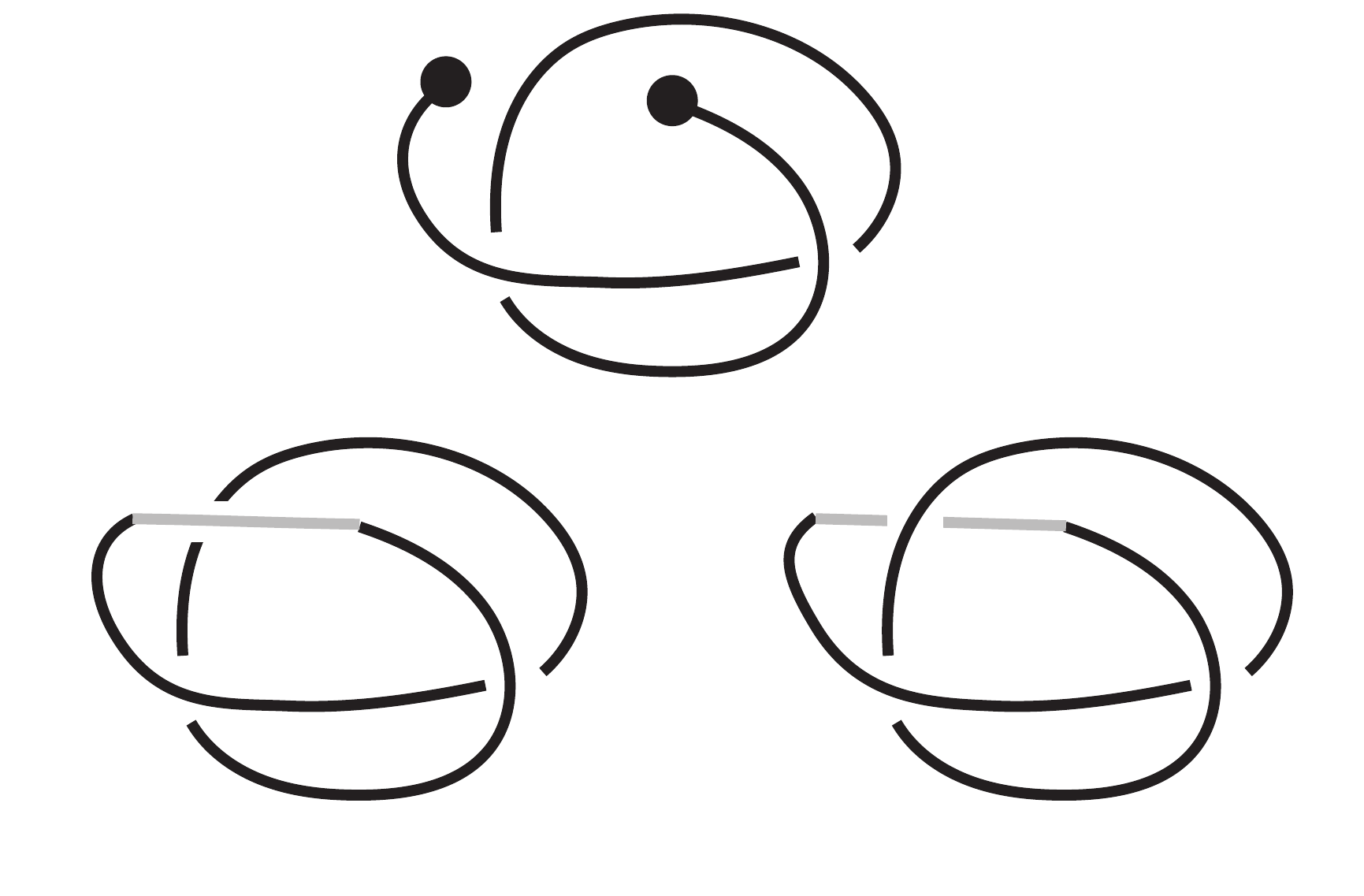}
\caption{The over and underclosure, respectively, of knotoid $2_1$. The overclosure gives the trivial knot; the underclosure gives the trefoil. The closure arc is shown in grey.}
\label{spherical_vs_planar}
\end{figure}
\end{center}

\subsection*{Acknowledgements} The research was supported by Williams College and NSF Grant DMS-1947438 supporting the SMALL Undergraduate Research Project.

\section{Spherical Knotoids} \label{spherical_knotoid_section}

\subsection{Maps from Spherical Knotoids to Knots}
We begin by describing two maps taking spherical knotoids to knots in the thickened torus $T \times (0, 1)$, the set of which we denote by $\mathcal{K}(T \times (0, 1))$. 

\begin{figure}[htbp]
\centering
\includegraphics[scale=0.4]{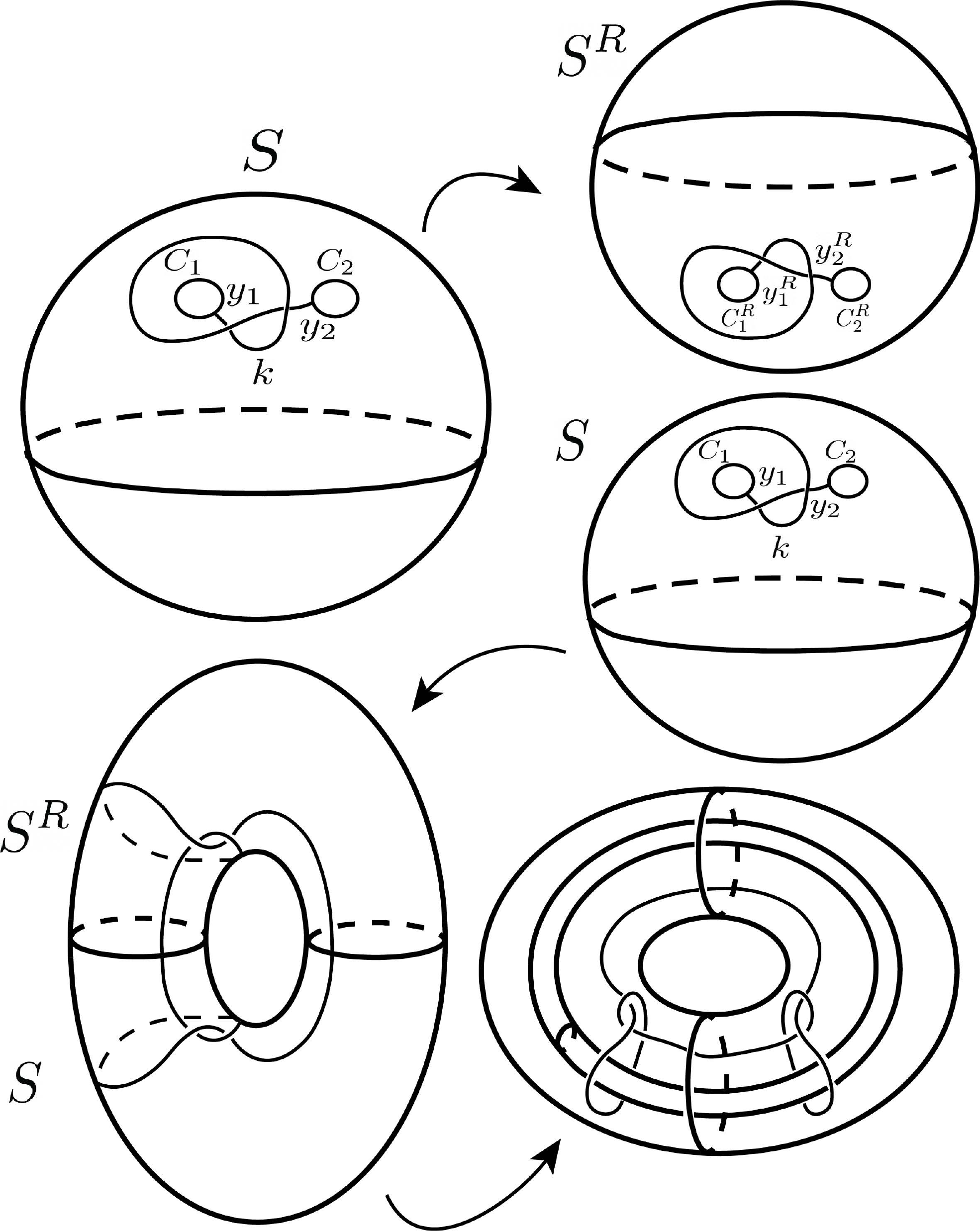}
\caption{The procedure described in Definition \ref{defnSphericalDoublingMap}.}
\label{sphericaldoubleFigure}
\end{figure}

\begin{definition} \label{defnSphericalDoublingMap}
Consider a knotoid diagram $\mathcal{D}$ on $S^2$ representing the spherical knotoid $k$. Denote the two endpoints of $k$ by $x_1$ and $x_2$, and by abuse of notation, denote the immersed arc by $k$. Take disk neighborhoods $D_i$ of the $x_i$ for $i = 1, 2$, chosen sufficiently small so that the circle boundaries $C_i := \partial D_i$ are punctured exactly once by $k$. Denote the punctures by $y_i$. Consider the surface $S= S^2 \setminus (\mathring{D_1} \cup \mathring{D_2})$ and take a reflected copy $S^R$, which contains reflected copies $C_i^R$ once-punctured by $y_i^R$. Now glue $S$ to $S^R$ via an orientation-preserving homeomorphism from $C_i$ to $C_i^R$ such that $y_i$ and $y_i^R$ are identified, for $i = 1, 2$.  This yields a knot diagram on a torus. The theory of knots in a thickened surface $\Sigma \times I$ is equivalent to the theory of knot diagrams on $\Sigma$. Then we obtain a knot $K$ in a thickened torus $T \times I$. Finally, remove the torus boundaries to obtain $K$ living in $T \times (0, 1)$. This defines a map $\phi_{S^2}^D: \mathbb{K}(S^2) \to \mathcal{K}(T \times (0, 1))$ which we call the \textit{spherical reflected doubling map} (on knotoids). We write $\phi_{S^2}^D(k) = K$ and refer to $K$ as the image of $k$ under the reflected doubling map. See Figure \ref{sphericaldoubleFigure}. 
\end{definition} 

\begin{figure}[htbp]
\includegraphics[scale=0.5]{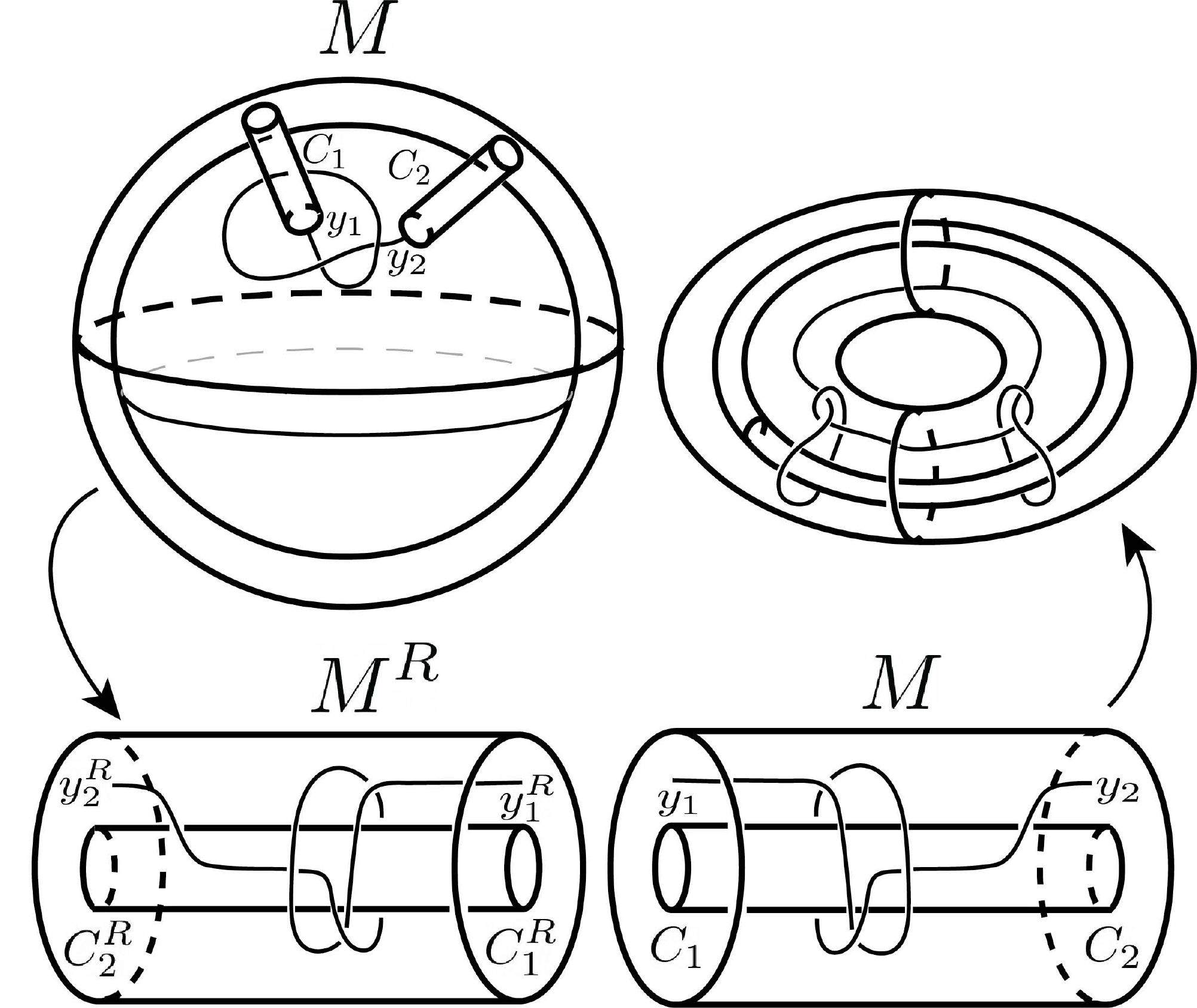}
\caption{An equivalent way of visualizing the construction in Definition \ref{defnSphericalDoublingMap}.}
\label{equivalentfigureDoubling}
\end{figure}

Equivalently, we may visualize this construction in the spherical rail diagram corresponding to $k$. That is, we can start with a thickened sphere $S^2 \times I$ together with two rails $\ell_i = \{x_i\} \times I$ for $i = 1, 2$, and an embedded arc $k$ with each endpoint in one of the rails. Let $M = (S^2 \times I) \setminus \mathring{N}(\ell_1 \cup \ell_2)$ be the manifold which is the thickened sphere with neighborhoods of the rails removed. Note that $M$ is just a thickened cylinder with a properly embedded arc. Let $C_i = \partial N(\ell_i)$ be the cylindrical boundaries of the rail neighborhoods, so that each contains one endpoint of $k$. Take a reflected copy $M^R$ and glue $M$ to $M^R$ along the $C_i$ via the orientation-preserving homeomorphism which is just a reflection, such that the endpoints from $k$ are identified. We might visualize this as pulling the $C_i$ up to the surface of $M$ before gluing in a ``natural way". See Figure \ref{equivalentfigureDoubling}. 

We remark that we remove the torus boundaries because we are interested in putting a hyperbolic structure on the corresponding knot complement in the thickened torus. This is not possible for a manifold which contains torus boundary components, as such tori prevent the existence of a hyperbolic metric. 

Now we may extend the notion of hyperbolicity via this construction for  knotoids. 

\begin{definition}
Let $k \in \mathbb{K}(S^2)$ be a spherical knotoid, and let $K = \phi_{S^2}^D(k)$ be the knot in the thickened torus which is the image of $k$ under the reflected doubling map. We say that $k$ is \textit{hyperbolic} (under the reflected doubling map) if $K$ is hyperbolic in the thickened torus, that is, $(T \times (0, 1)) \setminus \mathring{N}(K)$ admits a complete metric of constant sectional curvature -1. If $k$ is hyperbolic, then by the Mostow-Prasad Rigidity Theorem, we may define its \textit{hyperbolic volume} $\vol_{S^2}^D(k)$ (under the spherical reflected  doubling map) as
$$
\vol_{S^2}^D(k) := \frac{1}{2} \vol((T \times (0, 1)) \setminus \mathring{N}(K)). 
$$
\end{definition}



\begin{figure}[htbp]
\includegraphics[scale=0.4]{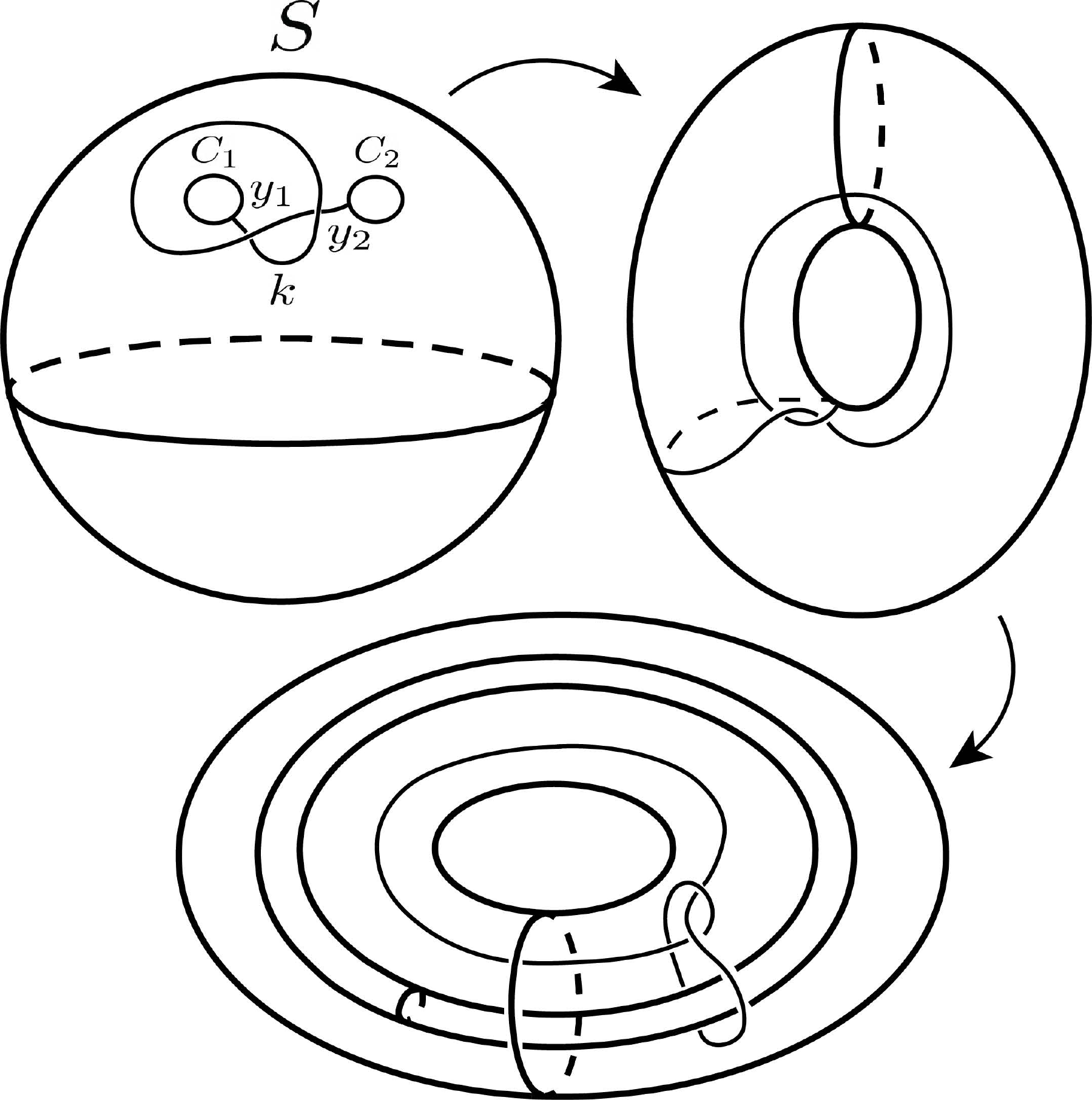}
\caption{The gluing procedure described in Definition \ref{defnSphericalGluingMap}.}
\label{sphericalgluingFigure1}
\end{figure}

Instead of doubling across the circle boundaries $C_i$ by taking a reflected copy of $S$ as in Definition \ref{defnSphericalDoublingMap}, we might instead choose to simply glue $C_1$ to $C_2$. This yields a second map from $\mathbb{K}(S^2)$ to $\mathcal{K}(T \times (0, 1))$ as in the following definition: 

\begin{definition} \label{defnSphericalGluingMap}
Consider a knotoid diagram $\mathcal{D}$ on $S^2$ representing the spherical knotoid $k$. We keep the notation $x_i$, $D_i$, $C_i$, $S$, and $y_i$ as in Definition \ref{defnSphericalDoublingMap}. Instead of doubling, simply glue $C_1$ to $C_2$ via an orientation-reversing homeomorphism such that $y_1$ and $y_2$ are identified. This yields a knot diagram on a torus which is equivalent to a knot $K$ in the thickened torus $T \times I$. Finally, remove the torus boundaries to obtain $K$ living in $T \times (0, 1)$. This defines a map $\phi_{S^2}^G: \mathbb{K}(S^2) \to \mathcal{K}(T \times (0, 1))$ which we call the \textit{spherical gluing map} (on knotoids). We write $\phi_{S^2}^G(k) = K$ and refer to $K$ as the image of $k$ under the gluing map. See Figure \ref{sphericalgluingFigure1}.
\end{definition}

\begin{figure}[htbp]
\centering
\includegraphics[scale=0.4]{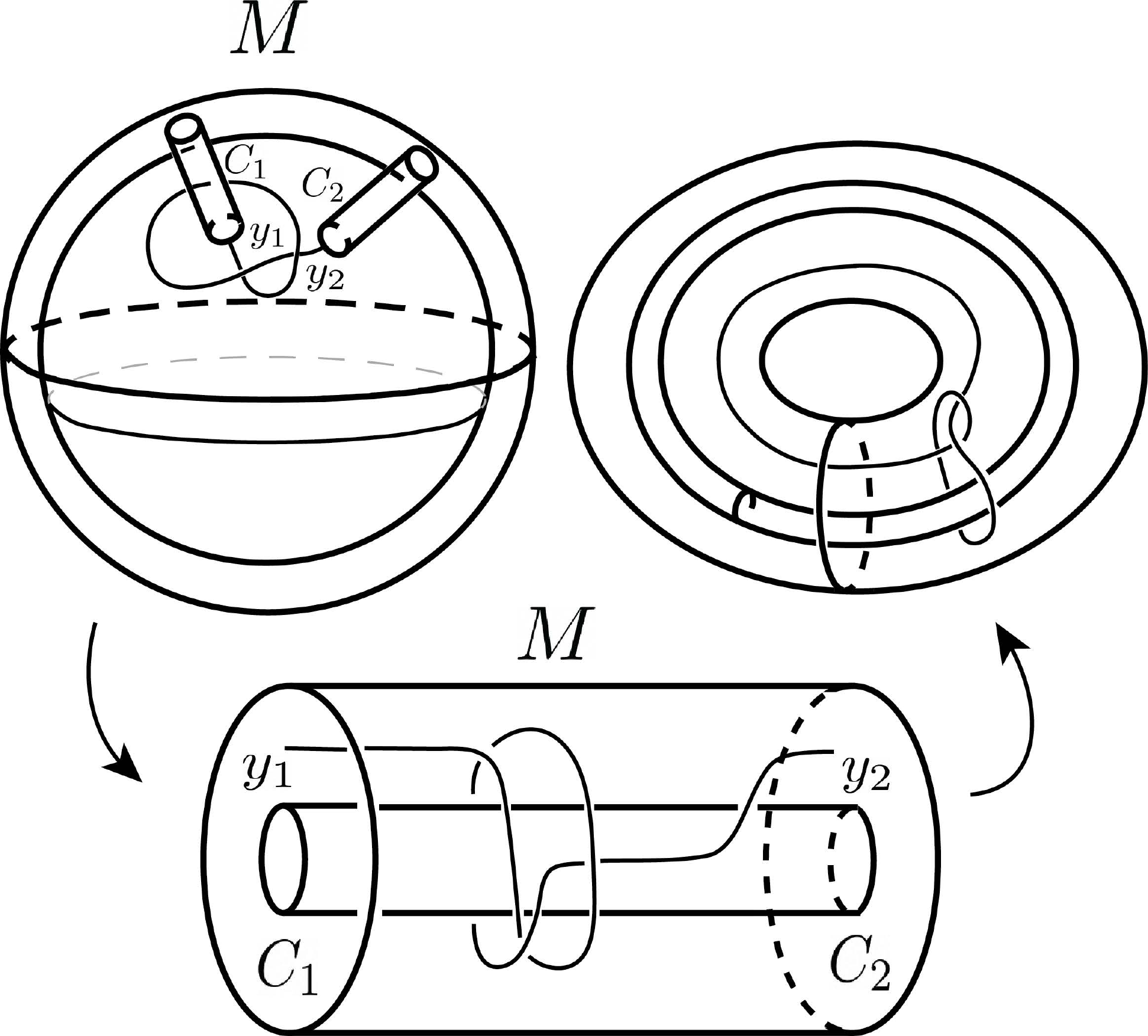}
\caption{An equivalent way to visualize the gluing in Definition 2.3}
\label{equivalentwaygluingFigure}
\end{figure}

This construction dates back to Turaev \cite{turaev} in his initial paper defining knotoids. Again, there is an equivalent way to visualize the gluing in the corresponding spherical rail diagram. As before, start with a thickened sphere $S^2 \times I$ with two rails $\ell_i = \{x_i\} \times I$ for $i = 1, 2$, and an embedded arc $k$. Let $M$ be $(S^2 \times I) \setminus \mathring{N}(\ell_1 \cup \ell_2)$, which is just a thickened cylinder with an embedded arc. Let $C_i = \partial N(\ell_i)$ as before, each containing one puncture from $k$. Now glue $C_1$ to $C_2$ via the ``natural" orientation-reversing homeomorphism such that the punctures are identified. See Figure \ref{equivalentwaygluingFigure}.

The same remark on removing the torus boundaries applies here. Also note that if $k$ is an oriented spherical knotoid, $\phi_{S^2}^G$ takes $k$ to an oriented knot in $T \times (0, 1)$. 

In \cite{siberian}, the authors show that the map $\phi_{S^2}^G(k)$ is well defined. It is also a main result of that paper that this map is injective for all prime knotoids of height greater than 1. 

In \cite{virtualclosure}, the authors point out that taking the virtual closure of a knotoid $k$ yields a knot in the thickened torus, and furthermore, this knot is isotopic in the thickened torus to the knot $K = \phi_{S^2}^G(k)$ obtained by the spherical gluing map. 

Now we may make another definition of hyperbolicity for spherical knotoids. 

\begin{definition} 
Let $k \in \mathbb{K}(S^2)$ be a spherical knotoid, and let $K = \phi_{S^2}^G(k)$ be the knot in the thickened torus which is the image of $k$ under the gluing map. We say $k$ is \textit{hyperbolic} (under the gluing map) if $K$ is hyperbolic in the thickened torus. If $k$ is hyperbolic, then its \textit{hyperbolic volume} $\vol_{S^2}^G(k)$ (under the spherical gluing map) is given by
$$
\vol_{S^2}^G(k) := \vol((T \times (0, 1)) \setminus \mathring{N}(K)). 
$$
\end{definition}

The second map yields similar results with regard to hyperbolicity as does the first.  

\begin{theorem}
Let $k \in \mathbb{K}(S^2)$ be a spherical knotoid, and consider its image under $\phi_{S^2}^G$ and $\phi_{S^2}^D$. Then $\phi_{S^2}^G(k)$ is hyperbolic if and only if $\phi_{S^2}^D(k)$ is hyperbolic. Furthermore, 
$$
\vol_{S^2}^D(k) = \vol_{S^2}^G(k). 
$$
\end{theorem}

\begin{proof}
Let $K^D$ be the knot in the thickened torus obtained from the reflected doubling map, and let $M^D = (T \times (0, 1)) \setminus N(K^D)$ be its complement in the thickened torus. Similarly, let $K^G$ be the knot in the thickened torus obtained from the gluing map, and let $M^G = (T \times (0, 1)) \setminus N(K^G)$ be its complement in the thickened torus. We claim that we can cut and reglue $M^D$ to obtain the disjoint union of $M^G$ and a reflected copy $(M^G)^R$. By construction, $M^D$ contains two once-punctured open annuli, $C_1$ and $C_2$ (as in Definition \ref{defnSphericalDoublingMap}), which we may think of as two thrice-punctured spheres since they don't contain their boundaries. After cutting along the $C_i$, we obtain a copy of $M$ and $M^R$, reusing notation from Definition \ref{defnSphericalDoublingMap}, which contain thrice-punctured spheres $C_1, C_1^R$ and $C_2, C_2^R$ respectively. Glue $C_1$ to $C_2$ and $C_1^R$ to $C_2^R$ as in Definition \ref{defnSphericalGluingMap}, yielding the disjoint union of $M^G$ and $(M^G)^R$. 

In \cite{thricepunctured}, the author proves that if a finite volume hyperbolic manifold $M$ can be cut and reglued along incompressible, embedded thrice-punctured spheres, resulting in a new manifold $M'$, then $M'$ is also finite volume hyperbolic and furthermore, $\vol(M) = \vol(M')$. Hence, if $M^D$ is hyperbolic, then $M^G$ is hyperbolic, and by reversing the cut-and-paste instructions, we get the converse. Equality of volumes also follows: 
\begin{align*}
    \vol_{S^2}^D(k) &= \frac{1}{2}\vol(M^D) \\
    &= \frac{1}{2}\vol(M^G \sqcup (M^G)^R) \\
    &= \vol(M^G) \\
    &= \vol_{S^2}^G(k). 
\end{align*}

\end{proof}



\begin{figure}[htbp]
\centering
\includegraphics[scale=0.4]{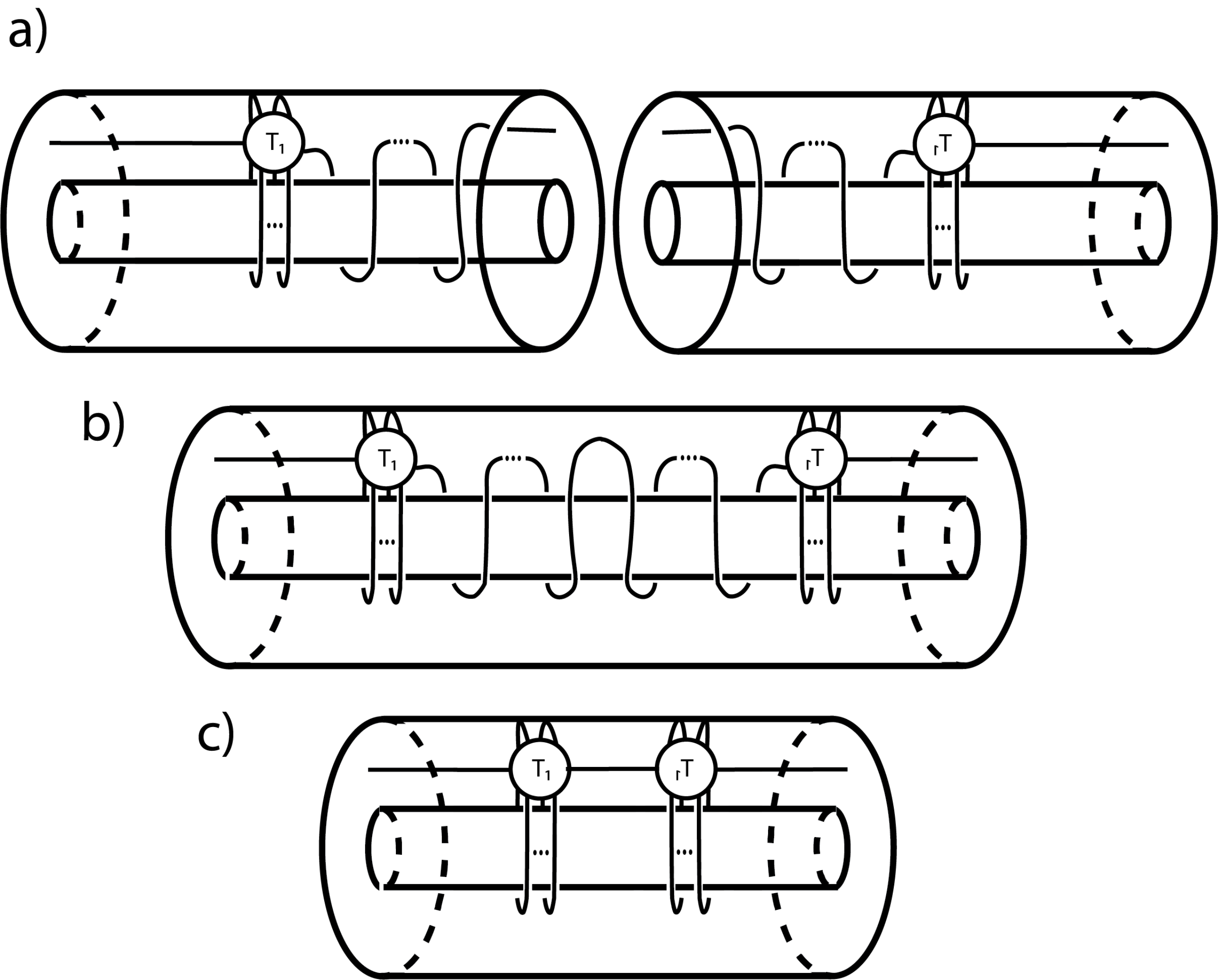}
\caption{The doubling construction undoes any twisting of the knot around the core.} 
\label{twists}
\end{figure}

Henceforth, we may say that $k$ is hyperbolic without referring to the reflected doubling or gluing maps, and we simply denote its volume by $\vol_{S^2}(k)$. The gluing map yields a simpler knot in the thickened torus and hence is easier to compute with. However, there are several advantages to the reflected doubling map. There is a distinction between isotopy of knots in the thickened torus versus homeomorphism of knot complements in the thickened torus.  For example, we might wind the embedded arc $k$  around the core curve of the thickened cylinder some number of times before gluing the boundaries of the cylinder together as in the spherical gluing map. The resulting knot would have complement in the thickened torus homeomorphic to the complement of the knot obtained without winding, but the two knots are not isotopic via ambient isotopy in the thickened torus. However, the reflected doubling action undoes any arbitrary twisting around the core (see Figure \ref{twists}).  

Another advantage to the reflected doubling map is it extends to a map on \textit{generalized knotoids}, which we define in \cite{generalizedknotoids}. Then we are able to define hyperbolicity for this larger class of objects. 

We note that a different definition of hyperbolicity for spherical knotoids was introduced in \cite{BBHL}. There, the initial knotoid was realized with a rail diagram in $S^2 \times I$.  Taking the branched double cover over the rails yields a knot in $S^2 \times I$, which after capping off the spherical boundaries, becomes a knot in $S^3$. Then if the resultant knot is hyperbolic, the volume associated to the knotoid is the volume of that knot. Note that the following proposition holds for both our construction and that construction.

\begin{proposition} \label{nonhyperbolicity_of_knottype} A knotoid that is knot-type is never hyperbolic, and a knotoid that is not knot-free is never hyperbolic.
\end{proposition}

\begin{proof}
By Thurston's Hyperbolization Theorem, $M$ is tg-hyperbolic if and only if $M$ has no essential spheres, disks, annuli, or tori. In the case of a knot-type knotoid, there is a disk in $S^2 \times I$ with boundary made up of four arcs, two on the rails and one each on $S^2 \times\{0\}$ and $S^2 \times \{1\}$  that is unpunctured by the knotoid. It becomes an essential annulus in the thickened torus. In the case of a  non-knot-free knotoid, there is a twice-punctured sphere that comes from the composition of the knot with the knotoid. It is mapped to an essential annulus in the knot complement in the thickened torus. Both of these annuli obstruct hyperbolicity.





\end{proof}

\subsection{Hyperbolicity}

It is well known that the composition of two classical knots can never be hyperbolic. However, in the case of oriented spherical knotoids, the situation is quite different: 

\begin{proposition} The product of two  oriented hyperbolic spherical knotoids  is always hyperbolic, and hyperbolic volumes of spherical knotoids are additive under multiplication. 
\label{product_is_additive}
\end{proposition} 

\begin{proof}
Let $k_1, k_2 \in \mathbb{K}(S^2)$ be oriented hyperbolic spherical knotoids; we prove that the image $\phi_{S^2}^G(k_1k_2)$ in the thickened torus is hyperbolic. Since the knotoids are on the sphere, for convenience of visualization we can always choose the tail of $k_1$ and the head of $k_2$ to be in their respective exterior regions, as mentioned in the introduction. Thus, viewing the thickened torus as the complement of the Hopf link, the image of the product of $k_1k_2$ under the spherical gluing map appears as in Figure \ref{productisbeltedsum}. 

\begin{figure}[htbp]
\centering
\includegraphics[scale=0.4]{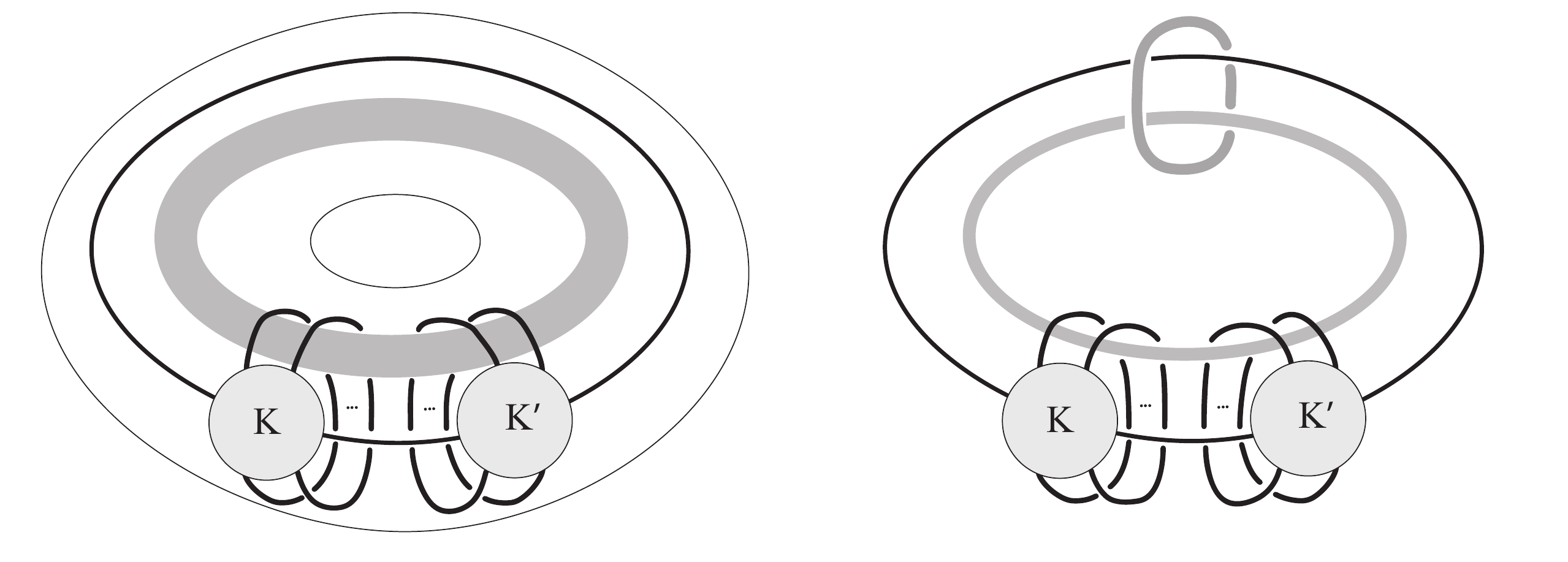}
\caption{On the left is the image of the composition of two knotoids in the thickened torus; on the right, we replace the thickened torus with the two Hopf Link components to make it clearer that the result is homeomorphic to a belted sum of two links.}
\label{productisbeltedsum}
\end{figure}

In \cite{thricepunctured},  the \emph{belted sum} $L_1 \#_b L_2$ of two links $L_1, L_2$, as shown in Figure \ref{beltedsum} is defined. There, it is proved that if $L_1$ and $L_2$ are both hyperbolic, then so is their belted sum, and  the volume of the belted sum is the sum of the volumes of $L_1$ and $L_2$. Figure \ref{productisbeltedsum} shows that when we replace the thickened torus with the complement of the Hopf link, $\phi^G_{S^2}(k_1k_2)$ is the belted sum of $\phi^G_{S^2}(k_1)$ and $\phi^G_{S^2}(k_2)$. Therefore, the result follows. 

\begin{center}
\begin{figure}[htbp]
\includegraphics[scale=0.4]{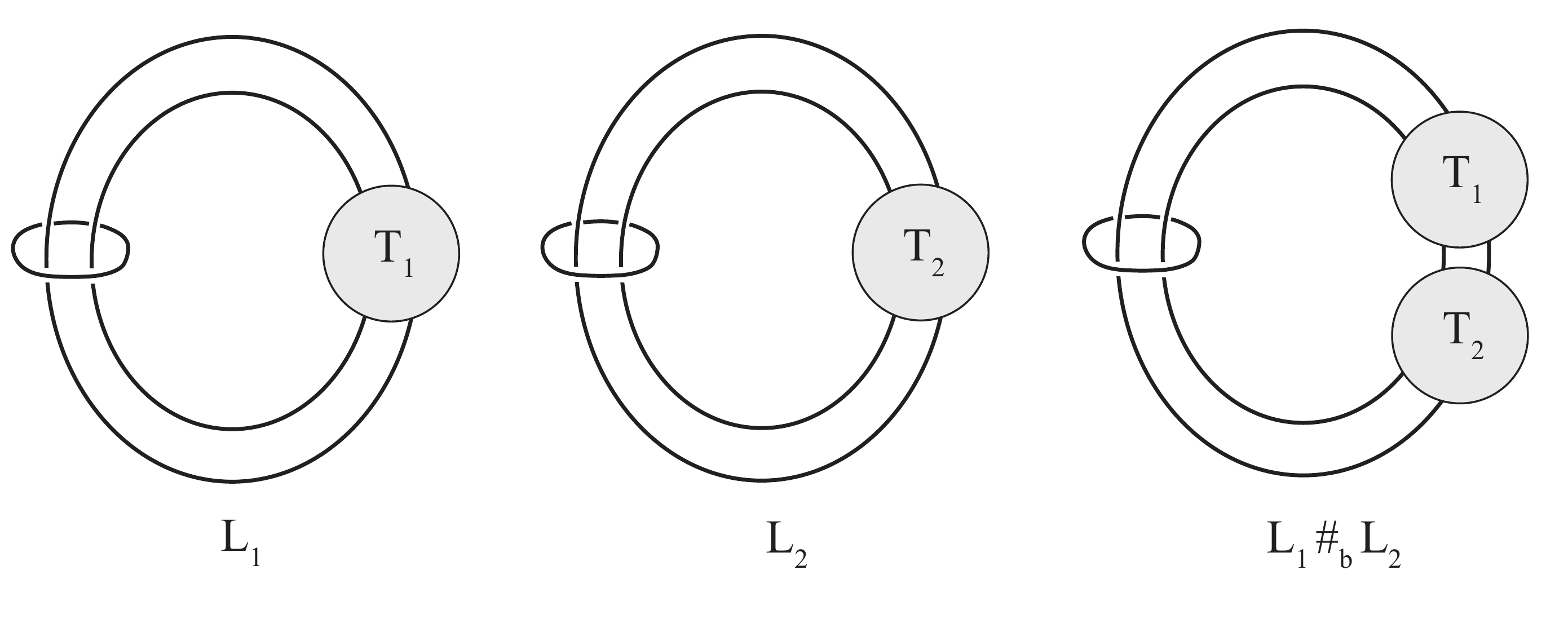}
\caption{The belted sum of two links $L_1$ and $L_2$.}
\label{beltedsum}
\end{figure}
\end{center}
\end{proof}

Note that in the case of the double branched cover construction from \cite{BBHL}, the product of two nontrivial spherical knotoids is never hyperbolic as there is always an essential annulus in the double branched cover.

We also find a large class of hyperbolic knotoids: 

\begin{definition} A knotoid $k$ (either in $S^2$ or $\mathbb{R}^2$) is \textit{closure-alternating} if, in a projection that realizes the minimum height, the crossings on a shortcut $s$ realizing the minimum height can be chosen so that $k \cup s$ yields a projection of an alternating knot. \label{alternating_knotoid_definition}
\end{definition}

\begin{definition}
We say a knotoid diagram in a projection surface $\Sigma$ is \textit{reduced} if there is no circle $\gamma$ which intersects the diagram exactly once through a double point corresponding to a crossing such that it is transverse to both strands at the intersection and $\gamma$ bounds a disk in $\Sigma$ not containing an endpoint of the knotoid.
\end{definition}

Note that given any knotoid diagram, we can always perform Reidemeister moves away from the endpoints to make it reduced. 

\begin{definition}
A knotoid $k$ in $S^2$ (resp. $\mathbb{R}^2$) represented by a reduced diagram is \textit{weakly knot-free} if for any circle $\gamma \subset S^2$ (resp. $\mathbb{R}^2$) which intersects the diagram exactly twice and bounds a disk $D$ not containing an endpoint of $k$, $D$ intersects the diagram of $k$ in a single arc with no crossings. 

Similarly, a knot $K$ in $S^2$ (equivalently $\mathbb{R}^2$) with a reduced diagram is \textit{weakly prime} if for any circle $\gamma \subset S^2$ which intersects the diagram exactly twice and bounds a disk $D$, $D$ intersects the diagram in a single arc with no crossings. 
\end{definition}

\begin{theorem} Let $k \in \mathbb{K}(S^2)$ be a weakly knot-free, closure-alternating spherical knotoid of height $1$. Then $k$ is hyperbolic. \label{Hyperbolicity_for_height_1_alternating_spherical}
\end{theorem}

\begin{proof} 
Choose a reduced closure-alternating projection of $k$. Let $\alpha$ be a shortcut which has no crossings with itself, connects the two endpoints, and intersects $k$ in one point. Let the knot $K \in \mathcal{K}(S^2)$ be $\alpha \cup k$. Furthermore, by our definition of closure-alternating, we can choose whether $\alpha$ crosses over or under the strand of $k$ it intersects such that the resulting diagram of $K$ is alternating. 

We show $K$ is weakly prime. If not, then there is a circle $C$ intersecting the projection at two points, with crossings to both sides. The two intersections cannot both be with $\alpha$, as the only crossing with $\alpha$ also involves $k$, so there would be a trivial sub-arc of $\alpha$ to one side and no crossings to that side of $C$. If both intersections are with $k$, then $\alpha$ lies entirely to one side of $C$.  There must be crossings of $k$ with itself to the  other side which contradicts the fact $k$ is weakly knot-free.   If one intersection is with $\alpha$ and the other with $k$, then to one side of $C$ is the single crossing between $\alpha$ and $k$ and to the other side $\alpha$  has no crossings. Hence, we can slide the intersection point with $C$ along $\alpha$ until we pass the endpoint and move the intersection point to be on $k$, contradicting our previous case.
Since $K$ is weakly prime, by \cite{menasco}, $K$ is prime.\\ 

We consider two cases separately, as to whether or not  $K$ is a (2,q)-torus link.  
In the case $K$ is not a $(2,q)$-torus knot, we use a theorem from \cite{colin!}.

Given a reduced connected alternating diagram $\mathcal{D}$ of a knot or link $K$, we can add a trivial component $J$ that intersects the projection plane in two points in two distinct non-adjacent regions of the projection plane such that disk it bounds is punctured twice by $K$. If $\mathcal{D}$ is weakly prime and not a projection of a 2-braid knot or link, then $L = K \cup J$ is hyperbolic.

\begin{figure}[htbp]
\centering
\includegraphics[scale=0.4]{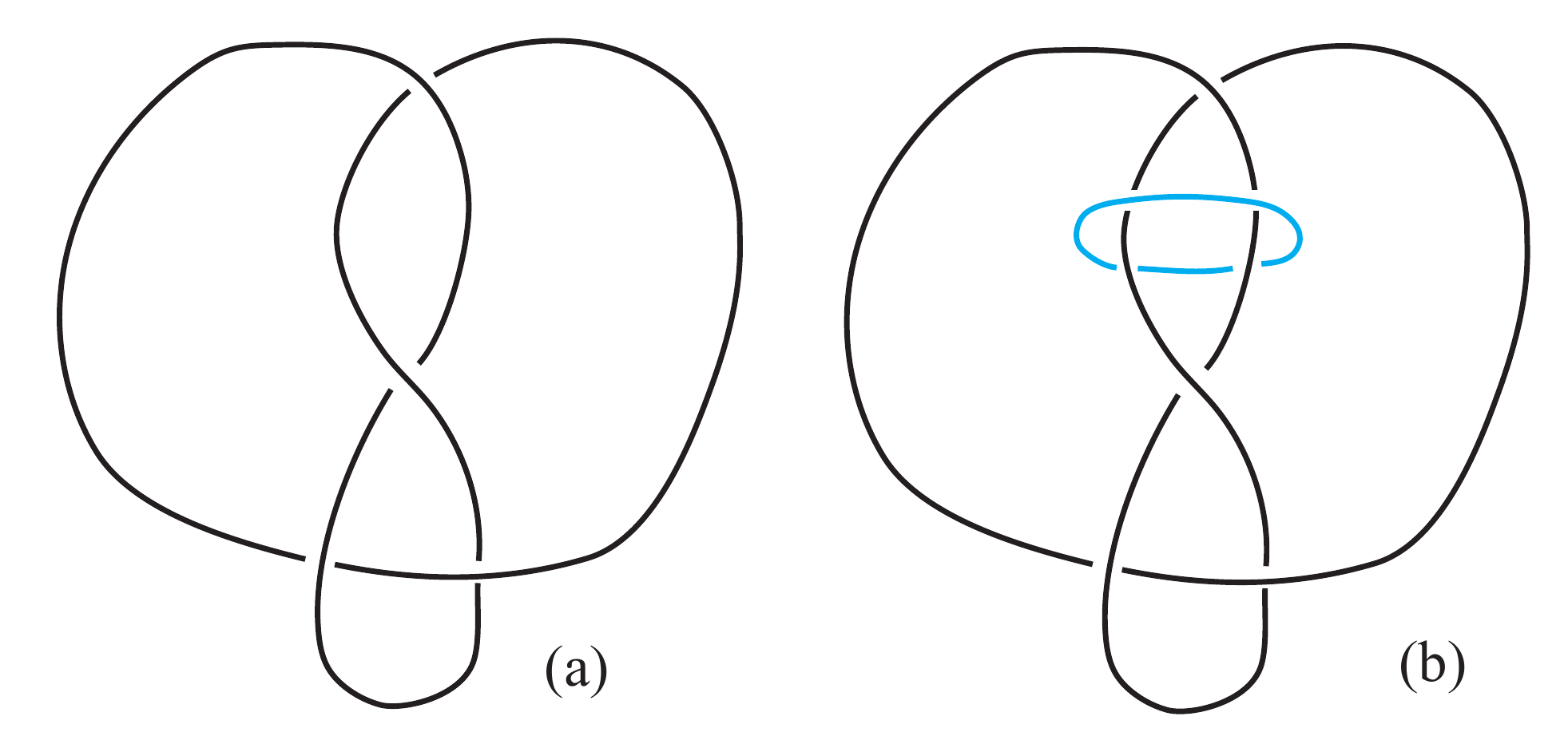}
\caption{To the left is the figure-eight knot, and to the right is an augmented alternating link obtained from the figure-eight knot by adding a trivial component.}
\label{augmentedlink}
\end{figure}



In our case, $K$ meets the hypotheses of the theorem,  so we augment $K$ with a component $J$ that passes through two of the regions meeting at the new crossing that are kittycorner to one another as in Figure \ref{usingchainmove} lower left to obtain a hyperbolic link. 

Now, in \cite{chainlemma}, the authors define a \textit{chain move} on a link $L$ as follows: start with a trivial component of $L$ bounding a disk twice punctured by $L$. As in Figure \ref{chainmove}, replace the tangle on the left with the tangle on the right to obtain a link $L'$. Chain moves in \cite{chainlemma} are defined more generally, but this specific construction is sufficient for our purposes.

\begin{figure}
    \centering
    \includegraphics[scale=0.5]{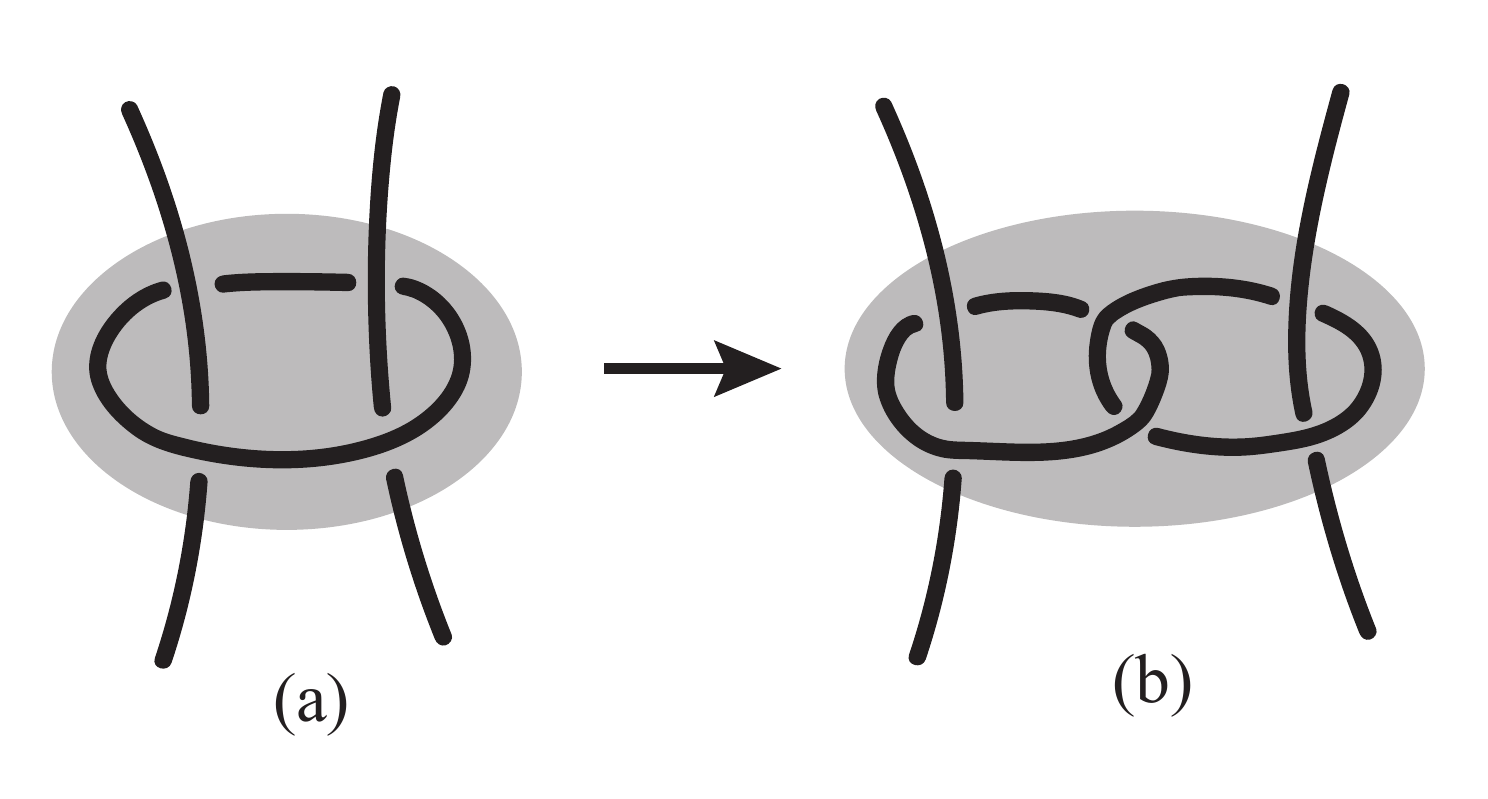}
    \caption{An example of a chain move on a link.}
    \label{chainmove}
\end{figure}

In Theorem 3.1 of \cite{chainlemma}, it is proved that if $L$ was originally hyperbolic in $S^3$, then so is $L'$, provided that if $\mathcal{B}$ is a ball surrounding the tangle as depicted in \ref{chainmove}, then $L \setminus (\mathcal{B} \cap L)$ is not the rational tangle $\infty$, 0, -1, or \\
-2 0 (this third being the vertical tangle -2).  However, our original knotoid must have at least two crossings to be the image of a single closed interval, so the closure must have at least three crossings, which prevents it from corresponding to any of these tangles. Therefore, after applying the move as depicted in Figure \ref{chainmove}, the resulting link $L''$ is also hyperbolic. See Figure \ref{usingchainmove}.

\begin{figure}
    \centering
    \includegraphics[scale=0.5]{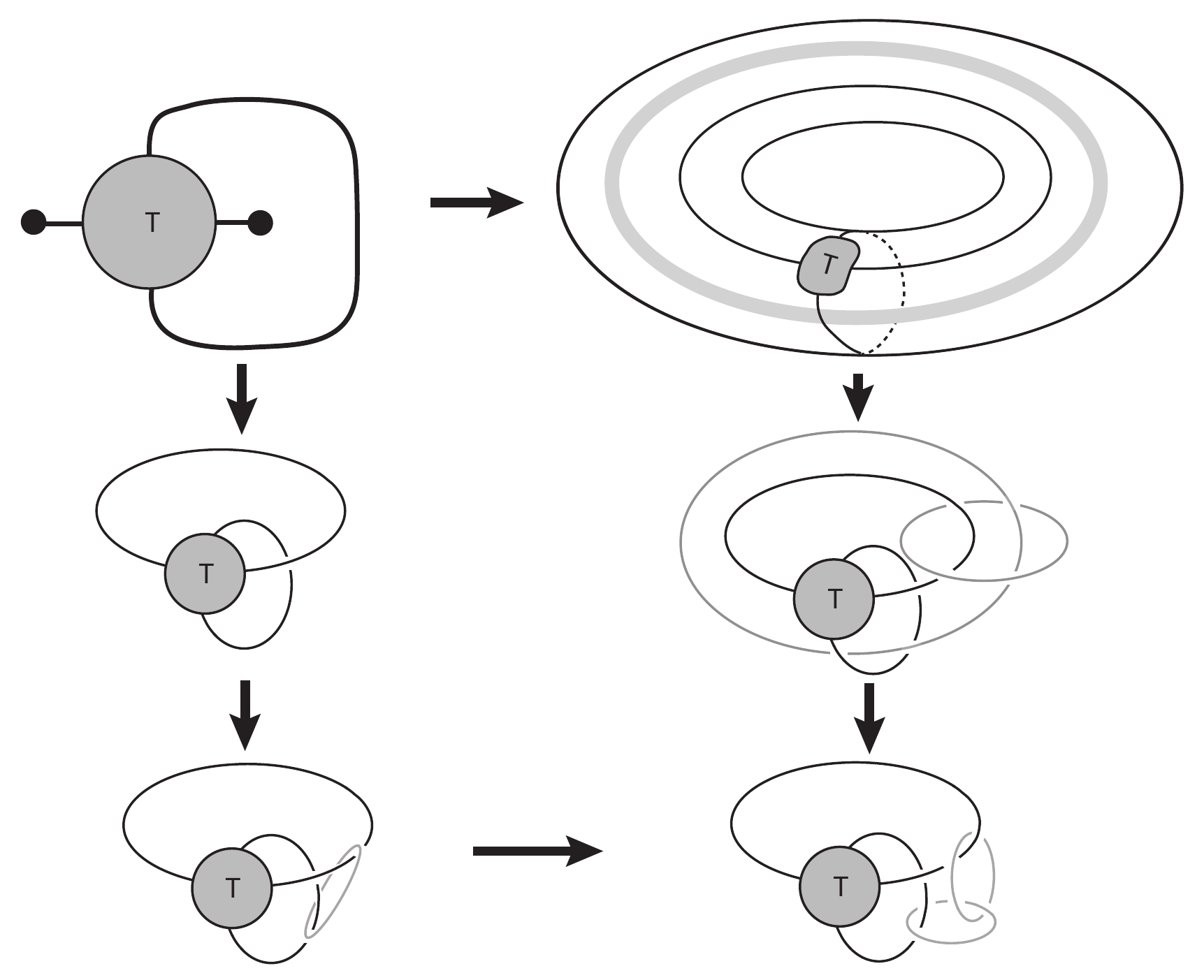}
    \caption{To obtain the knot in the thickened torus, which can be represented by the link in the lower right as we see by following the arrows left and then down in the figure, we can instead augment and then apply the chain move to get there as we follow the arrows down and then right in the figure.}
    \label{usingchainmove}
\end{figure}

Using once more the fact that the complement of the Hopf link in $S^3$ is homeomorphic to the thickened torus, the complement of $L''$ is homeomorphic to the complement of $\phi^D_{S^2}(k)$ in $T^2 \times I$. Since $L''$ is hyperbolic, so is $\phi^D_{S^2}(k)$, as desired. 

The other case to consider is when $K$ is a $(2,q)$-torus knot. This corresponds to an integer knotoid $k$ with $|k| >2$ as defined below. The technique used to prove hyperbolicity applied above does not work in this case, as augmenting a $(2,q)$-torus knot need not result in a hyperbolic link. But we consider the link obtained by  adding in the two components corresponding to the Hopf link. That this new link is hyperbolic follows from Lemma 3.4 of \cite{chainlemma}. There, it is proved that a rational tangle glued to the tangle appearing in Figure \ref{chainmove}(b) is hyperbolic unless the rational tangle is $\infty$, 0, -1 or -2 0. In our case the tangle is an integer tangle $k$ with $|k| > 2$. (Note that we count the additional crossing coming from the shortcut.) Hence, it cannot be any of these four tangles. 
\end{proof}

\begin{figure}[htbp]
\centering
\includegraphics[scale=0.5]{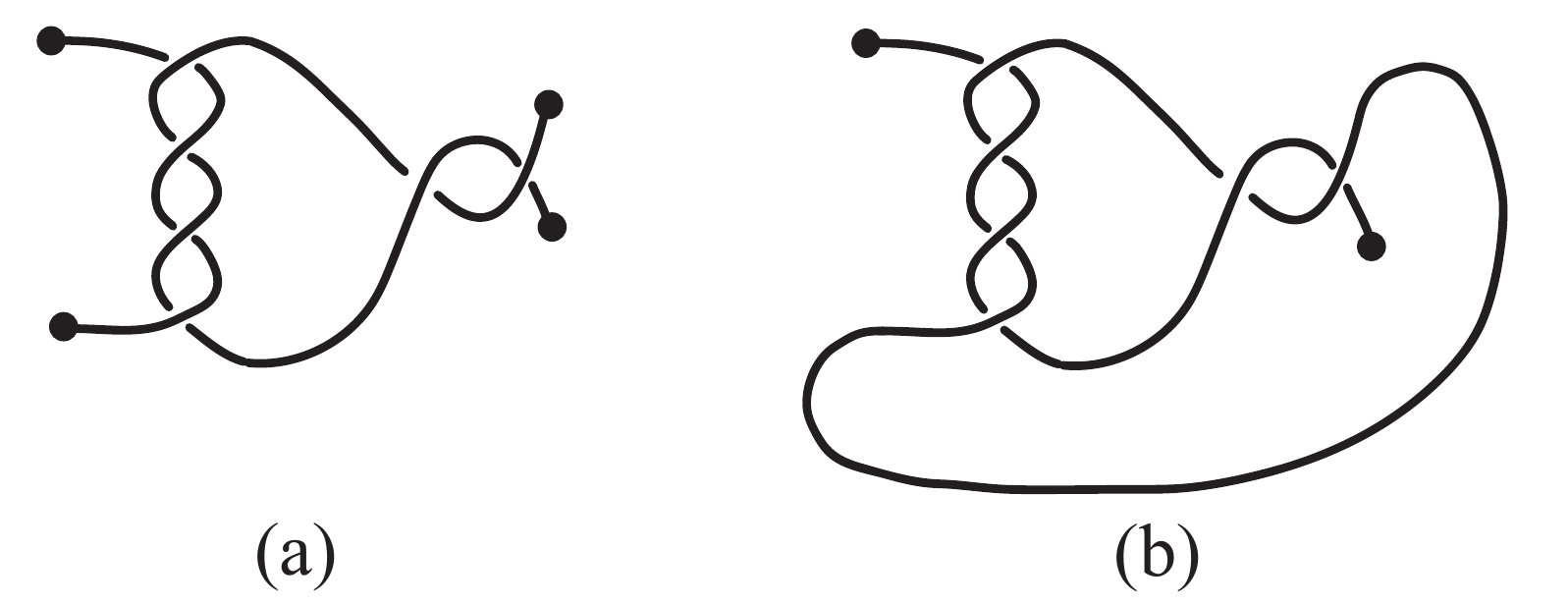}
\caption{To the left is a rational tangle, and to the right is  the rational knotoid we obtain from the tangle.}
\label{rationalknotoid}
\end{figure}

\begin{definition} A \emph{rational knotoid} is a knotoid obtained from a rational tangle as for instance appears in Figure \ref{rationalknotoid}(a) by gluing the northeast strand to the southwest strand so that no more crossings are created in the diagram as in Figure \ref{rationalknotoid}(b). Note that for knotoids on the sphere, it does not matter whether the strand connecting the two endpoints is to the upper left or  the lower right of the tangle; otherwise, by convention we join these two strands on the lower right of the tangle.
\label{Rational_Knotoid_definition}
\end{definition}

Note that the rational tangle must be a tangle such that the gluing of the northeast and southwest strands yields a single component rather than two components.

For example, $2_1$ is the rational knotoid with one twist region and two crossings in that region (where a \emph{twist region} in a knotoid diagram is a maximal chain of bigons, or a single crossing if no bigon touches that crossing). Note that a rational knotoid is always height 1.

\begin{definition}
An integer knotoid is a rational knotoid arising from a rational tangle with one twist region.\label{Integer_knotoid_definition}
\end{definition}

We immediately have the following corollary to Theorem \ref{Hyperbolicity_for_height_1_alternating_spherical}:

\begin{corollary}
All rational knotoids are hyperbolic.
\label{Hyperbolicity_of_rational_knotoids_corollary}
\end{corollary}

\begin{proof} 
Every rational tangle has a projection which is alternating, as shown, for instance,  in \cite{alternatingrational} and is knot-free. The corresponding rational knotoid is height 1 and closure-alternating so this follows immediately. 
\end{proof}

Given our constructions, we can make the following observation: 

\begin{proposition}\label{virtualclosureproposition}
For a knotoid $k$, if the overclosure and virtual closure of the knotoid are both hyperbolic, then the volume of the virtual closure is strictly larger than the volume of the overclosure. The same result holds replacing overclosure with underclosure.
\label{Volume_of_virtual_closure_and_overclosure}
\end{proposition}

\begin{proof}
In \cite{virtualclosure}, the authors point out that the map $\phi^G_{S^2}$ gives the same knot in the thickened torus as the virtual closure map. Using the fact that the complement of the Hopf link is homeomorphic to the thickened torus, we can describe our map as follows: start with a diagram of a knotoid $k$. Add an overstrand (resp. understrand) $\ell$ joining the endpoints of $k$. Then, to obtain the knot in the thickened torus, we add in two trivial components that are a Hopf link intertwined with this knot as depicted for example in Figure \ref{usingchainmove}. 

In other words, we take the over or under closure of $k$ and then add two additional components. But when the resulting link is hyperbolic, adding the extra components always strictly increases volume. Therefore, the claim follows.
\end{proof}


We note the following corollary which is unrelated to knotoids but of use in virtual knot theory.

\begin{corollary}
If $k$ is a classical knot, and both $k$ and the virtual knot $k'$ obtained by switching one crossing to virtual yields a non-classical virtual knot and both $k$ and $k'$ are hyperbolic, then the volume of $k$ is strictly less than the volume of $k'$. 
\label{Switching_crossing_to_virtual_increases_volume_corollary}
\end{corollary}

\begin{proof}

Note that the volume of a virtual link is discussed in \cite{SMALLvirtualknots}. Let $c$ be a crossing of a classical knot $k$. Cut the overstrand (resp. understrand) of $c$ right before and after $c$ to obtain a height 1 knotoid $K$. The overclosure (resp. underclosure) of $K$ is $k$, and the virtual closure is $k$ with $c$ switched to a virtual crossing. By Proposition \ref{virtualclosureproposition}, the volume of the latter is greater than the volume of $k$. 

\end{proof}


\section{Volume Bounds for Knotoids in $S^2$} \label{vol_bounds_spherical}

Explcit volumes in this section were all computed using the software SnapPy \cite{SnapPy}. 

It is a classical result of Thurston and J\o{}rgensen that volumes of hyperbolic 3-manifolds are well ordered in $\mathbb{R}$ (see, for instance, \cite{thurstonlecturenotes}). Thus, it makes sense to ask for the minimum volume of a knot, for example, or more generally an $n$-cusped link. In our case, it is relevant to ask for the minimum volume of a spherical knotoid. 

\begin{conjecture}\label{Knotoid_of_smallest_volume_conjecture} The unique spherical knotoid of smallest volume is the knotoid $2_1$. Its volume is approximately $5.33349\dots$ 
\end{conjecture}

Our conjecture is a weaker version of the conjectured minimum volume of a 3-cusped hyperbolic manifold (see, for instance, \cite{Agol2cusp}). Since a knot in the thickened torus may be thought of as a 3-component link where two components make up the Hopf link as described in the previous section, and since our conjectured minimum volume knotoid achieves the conjectured minimum volume of a 3-cusped hyperbolic manifold, $5.33349\dots$, proving this stronger conjecture would also prove our conjecture.

\begin{theorem}\label{Rational_knotoid_of_smallest_volume_lemma} The unique rational knotoid of smallest volume is the knotoid $2_1$, with volume approximately $5.33349\dots$. 
\end{theorem}

\begin{proof}

Throughout the proof, for convenience, we refer to the twist regions of a rational tangle with the ordering shown in Figure \ref{rationaltanglenotation} or its reflection about a NW-SE diagonal line.

\begin{figure}[htbp]
    \centering
    \includegraphics[scale=0.3]{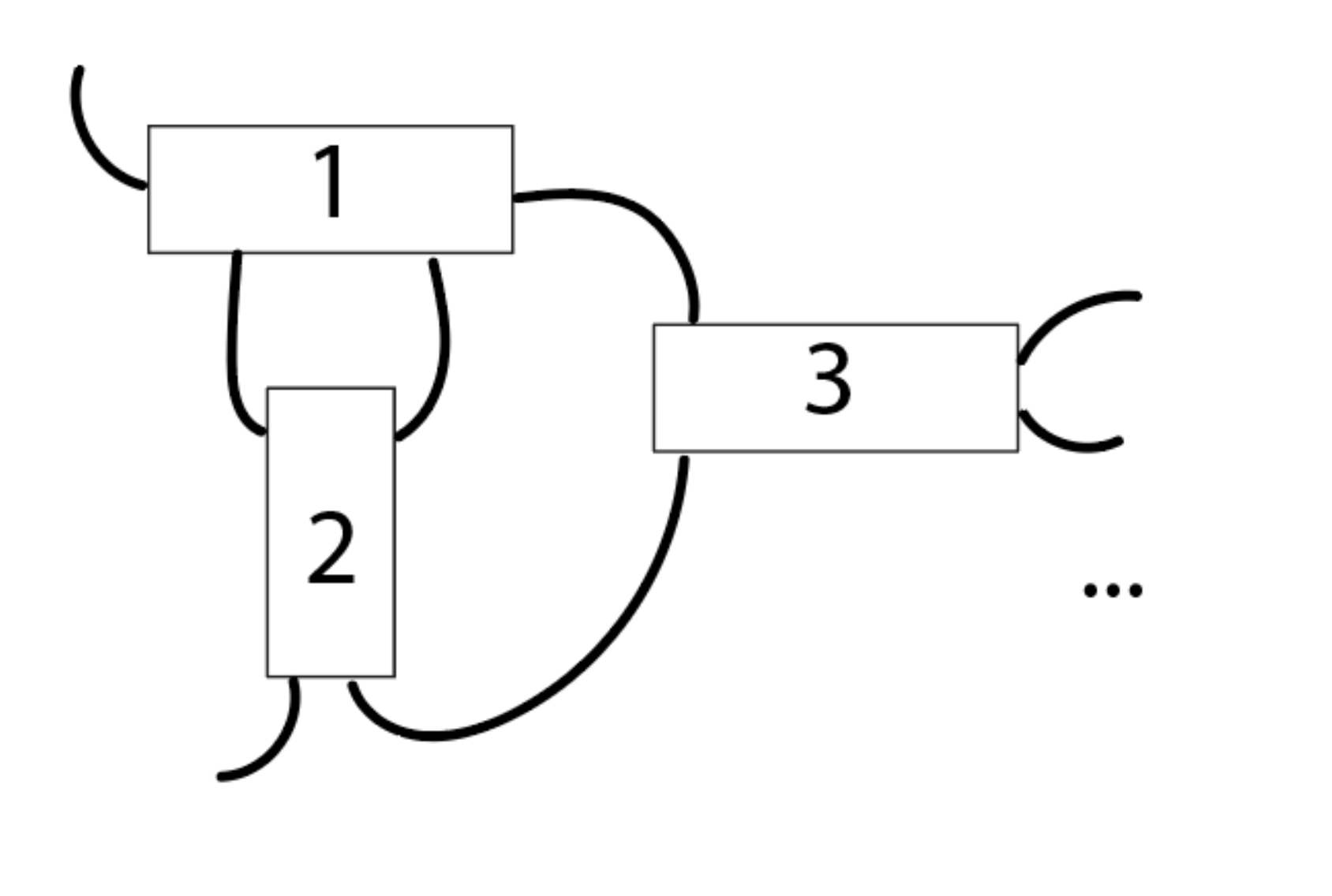}
    \caption{Our convention for ordering the twist regions of a rational tangle (the same convention is used for rational knotoids and as there, we allow for reflection across a NW-SE diagonal line.)}
    \label{rationaltanglenotation}
\end{figure}

We use again the fact that the complement of the Hopf link in $S^3$ is homeomorphic to the thickened torus. Therefore, computing the volume of a rational knotoid is equivalent to starting with an alternating projection of a rational tangle and gluing the northeast and southwest, and northwest and southeast, endpoints; these two new strands will intersect in one crossing, and we choose that crossing so that the diagram is still alternating. Finally, we add in the two Hopf link components representing the thickened torus as depicted in in Figure \ref{rationalinthickenedtorus}.

\begin{figure}[htbp]
    \centering
    \includegraphics[scale=0.25]{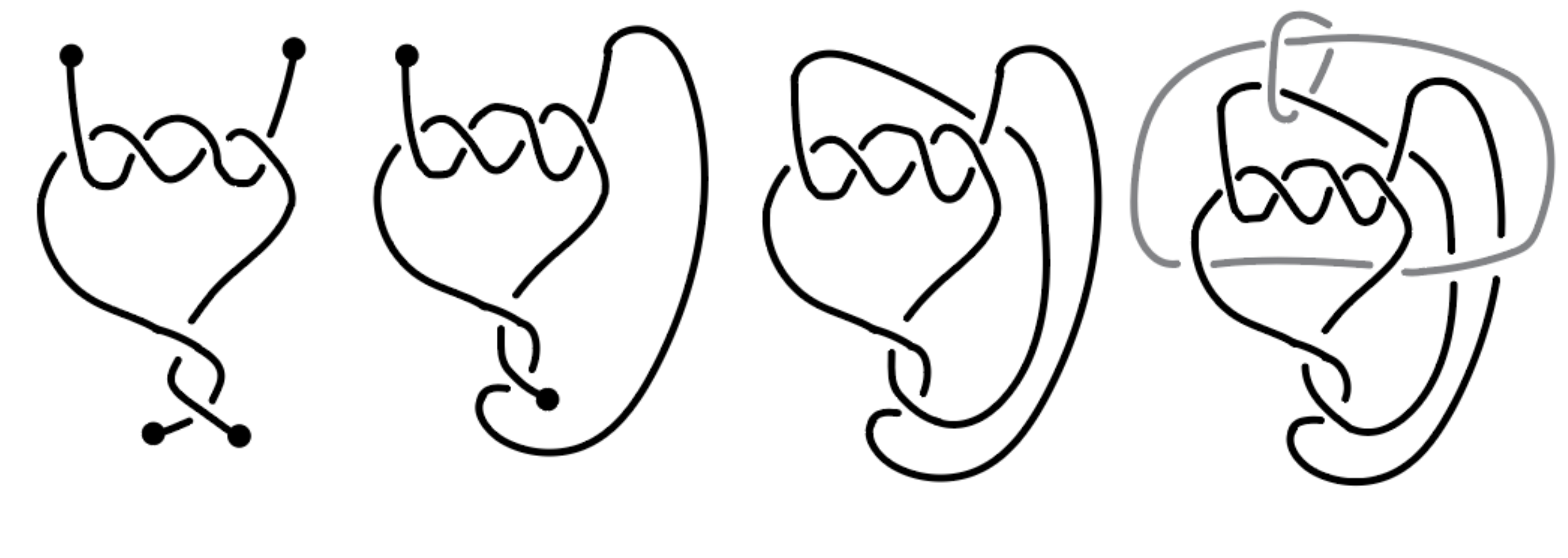}
    \caption{Similarly to the reasoning in Proposition 2.13, we can view the image of a knotoid in $\phi^D_{S^2}$ as a three component link. From left to right is a rational tangle, the associated rational knotoid, the image of the knotoid in the thickened torus, and the image of the knotoid in the thickened torus with two Hopf link components added to represent the thickened torus.}
    \label{rationalinthickenedtorus}
\end{figure}

Before the Hopf link components are added, the diagram is an alternating 2-bridge knot in $S^3$ (since the crossing we add when we join the endpoints of the tangle just adjoins onto the last twist region). By a result of Futer and Gu\'eritard from \cite{FuterGueritaud}, its volume is bounded below by

\begin{equation}\label{twobridgebound}
\begin{centering}
\operatorname{vol}(S^3- K) > 2v_{\operatorname{tet}}\operatorname{tw}(K) -2.7066,
\end{centering}
\end{equation}

where $v_{\operatorname{tet}} \approx 1.01494 \dots$ is the maximal volume of an ideal tetrahedron in hyperbolic 3-space and $\operatorname{tw}(K)$ denotes the number of twist regions in $K$. When we add in the two Hopf link components, the volume must increase, so this bound is also a lower bound for the volume of our knotoid. When there are at least 4 twist regions, plugging into this bound gives a lower bound on volume of approximately $5.4129\dots$ which already exceeds our claimed minimum, so we may restrict to the cases where there are 1, 2, or 3 twist regions. 

\subsection{The Case of $\operatorname{tw}(K) = 1$}

We first prove that all integer tangles other than the tangle with one twist region containing two crossings (which is the rational tangle whose corresponding rational knotoid is $2_1$) yield knotoids with volume greater than $5.33349\dots$. Note that if there is an odd number of crossings in the twist region, then according to our method of closing rational tangles as depicted in Figure \ref{rationalinthickenedtorus}, we do not get a knot in the thickened torus, but instead a 2-component link in the thickened torus.  Thus, we can assume that our twist region had an even number of twists. Given an integer tangle, we start by augmenting its twist region as depicted in Figure \ref{augmentedinteger}. 

\begin{figure}[htbp]
    \centering
    \includegraphics[scale=0.4]{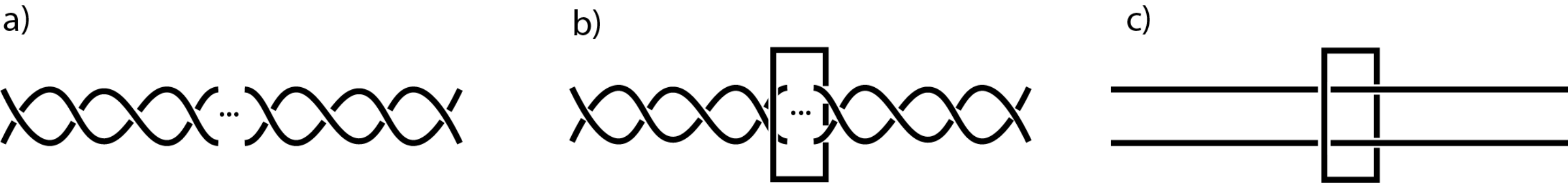}
    \caption{Augmenting a twist region and then removing crossings.}
    \label{augmentedinteger}
\end{figure}

By the results of \cite{Adamsetal} or \cite{FALsintorus}, the resulting link in the thickened torus will be hyperbolic.   

Since there is an even number of crossings in our twist region, we can cut along the twice-punctured disk bounded by $\ell$, untwist all of the crossings, and glue back together; this operation yields a homeomorphism of the link complement, preserving hyperbolicity and hyperbolic volume. We can use SnapPy to compute the volume of the resulting link $L$ in the thickened torus to be 10 $v_{tet} = 10.14942\dots$.

Because a reflection in the projection torus provides an isometry of the manifold preserving the longitude of the original cusp and reversing its meridian, a fundamental domain for the original cusp given by the meridian and longitude is a rectangle. 

Now, to obtain the original integer knotoid, we apply a Dehn filling to the augmented component. In \cite{dehnfillingbound}, the authors prove a lower bound on volumes for Dehn fillings: if $M$ is a hyperbolic link complement, with $C_1,\dots, C_k$ disjoint horoball neighborhoods of some subset of the cusps, let $s_1,\dots, s_k$ be the slopes of $\partial C_1,\dots, \partial C_k$ corresponding to our choice of Dehn fillings on the cusps, such that each length is greater than $2\pi$. Let $\ell_{\min}$ denote the minimum of these slopes. Then, the following bound holds: 

\begin{equation}\label{Dehnfillingboundeqn}
\begin{centering}
\operatorname{vol}(M(s_1,\dots, s_k)) \geq \left(1-\left(\frac{2\pi}{\ell_{\operatorname{min}}}\right)^2\right)^{\frac{3}{2}}\operatorname{vol}(M),
\end{centering}
\end{equation}

where $M(s_1, \dots, s_k)$ denotes the manifold resulting from Dehn filling along slopes $s_1,\dots, s_k$.

We see that the fundamental domain of the cusp neighborhood of the augmented component $\ell$ has longitude 4 and meridian $\sqrt{12}= 3.4641\dots$ when the cusp is maximally expanded. When we apply $(1, n)$ Dehn surgery, we can compute the corresponding slope by imagining gluing together $n$ copies of the fundamental domain side by side along their meridians and taking the length of the diagonal in the resulting $\sqrt{12}$ by $4n$ rectangle. The bound from Inequality \ref{Dehnfillingboundeqn} yields

\begin{equation}\label{onetwistslopebound}
\begin{centering}
\operatorname{vol}(M(s_1, \dots, s_n)) \geq \left(1-\left(\frac{2\pi}{\sqrt{12+16n^2}}\right)^2\right)^{\frac{3}{2}}10.149416.
\end{centering}
\end{equation}

(This bound is valid only when the value of $\sqrt{12+16n^2}$ exceeds $2\pi$, but this occurs for all $n \geq 2$). This volume bound exceeds $5.33349\dots$ for $n \geq 3$, so we only need to compare the $n=1$ case (which corresponds to $2_1$) and the $n=2$ case. However, we can directly compute that the volume for the $n=2$ case is greater than $5.33349\dots$ (SnapPy \cite{SnapPy} gives approximately $8.35550\dots$), so this case is complete. 

\subsection{The Case of $\operatorname{tw}(K) = 2$}

We follow the same idea as the first case, this time with  more casework. We start by augmenting each of the two twist regions of our rational tangle, cutting along the twice punctured disks bounded by our augmenting components, undoing the twists of each region, and gluing back together along these disks in the same way as the previous section. Unlike the $\operatorname{tw}(K) =1$ case, there are now several subcases: we can have no crossings left in either twist region, one crossing left in one twist region but none in the other, or a crossing left in each twist region. All of these resulting links have the same volume by \cite{colin!}, so for convenience, we assume that we are in the first case. 

We compute that the volume of the resulting link in the thickened torus is hyperbolic with volume $17.47714\dots$. It is possible to expand the cusps of the two augmenting components so that the dimensions of the fundamental domains are 4 by 2, and $\sqrt{3}+1 = 2.73205\ldots$ by 2.

Assume a $(1, n)$ Dehn filling along the first slope and a $(1, m)$ Dehn filling along the second slope. By the bound from Inequality \ref{Dehnfillingboundeqn}, the volume is lower-bounded by the maximum of the following:

\begin{equation}\label{Dehnfillingboundeqn2}
\begin{centering}
\operatorname{vol}(M) \geq \left(1- \left(\frac{2\pi}{(\sqrt{4+(4n)^2)}) }\right)^2\right)^{\frac{3}{2}}17.47714
\end{centering}
\end{equation}
and

\begin{equation}\label{Dehnfillingboundeqn3}
\begin{centering}
\operatorname{vol}(M) \geq \left(1- \left(\frac{2\pi}{\sqrt{4+((\sqrt{3}+1)m)^2}}^2\right)\right)^{\frac{3}{2}}17.47714...
\end{centering}
\end{equation}

\vspace{.15 in}

When $n > 2$ and $m > 3$, the volume exceeds $5.33349\dots$.

This subdivides our problem into several subcases: we can have 1, 2, or 3 crossings in the second twist region and arbitrarily many in the first, or we can have 1 or 2 crossings in the first twist region and arbitrarily many in the second. 

We describe in detail the computations for one subcase only (there is one twist in the first region and arbitrarily many in the second) to avoid being overly repetitive; the methods for the other subcases are identical. Augment the second twist region but not the first. The resulting link is hyperbolic with volume in the thickened torus $13.81328\dots$, and the cusp of the augmenting circle (which is no longer right angled) will have longitude $4$ and meridian $1 + (2\sqrt{3}+1)i$ when it is maximally expanded. Plugging these numbers into Inequalities \ref{Dehnfillingboundeqn2} and \ref{Dehnfillingboundeqn3}, we obtain a volume greater than $5.33349\dots$ when the number of twists is greater than 2. So, we restrict the possible knots with volume smaller than $5.33349\dots$ in this subcase to a set with two elements. We can compute by SnapPy that both volumes are greater than $5.33349\dots$ Therefore, this subcase is done, and the other subcases can be handled similarly.

\subsection{The Case of $\operatorname{tw}(K) = 3$}

This time there will be three levels of subcases, one in which there are three augmenting circles, one in which there are two, and one in which there is one; since the methods are the same as the previous two cases, to avoid being overly repetitive, we do not give as many details as before. (Furthermore, similarly to the case of $\operatorname{tw}(K) = 2$, there are additional cases depending on the parity of the number of twists in each region). We start with an augmented link with three crossing circles, which has volume approximately $24.80487\dots$. We can expand the cusp diagrams of the three crossing circles to have respective dimensions 4 by 2, 2 by 2, and 4 by $2\sqrt{3}+2$. Plugging into Equation \ref{Dehnfillingboundeqn2}, we see that if there is more than 1 twist in the first region, more than 3 twists in the second region, \emph{and} more than 1 twist in the third region, then the volume exceeds $5.33349\dots$.
\newline
This breaks our problem up into five more subcases: 

$$(1, a, b), (a, b, 1), (a, 1, b), (a, 2, b), (a, 3, b)$$

Each of these cases is handled by adding in two crossing circles to the two twist regions which can have arbitrarily many twists and using the usual Dehn filling bound from there. Each of the five subcases will break up into several more subcases: one set of subcases where we augment the twist region having $a$ twists and choose a constant number of twists in the twist region having $b$ twists, and one where we augment the twist region having $b$ twists and choose a fixed number of twists in the twist region having $a$ twists. Finally, as before, each of these remaining subcases is eventually restricted to a finite set of possible knotoids. We may check by hand that each volume is actually greater than $5.33349\dots$, and we are done. 
\end{proof}

\begin{question}

Are the volumes of all integer tangles less than the volumes of all other rational tangles? Note that this is equivalent to asking whether the volume of every non-integer rational tangle is greater than $10.149\dots$ The same methods as above would theoretically solve this problem, but much more casework and computational power is required.

\end{question}

\begin{question}

Is the minimum volume of an alternating knotoid $5.33349\dots$? This problem is also in theory solvable using the methods above, but the required casework is outside the scope of this paper. 

\end{question}

\begin{question} \label{Do_volumes_of_integer_tangles_increase_monotonically_problem}
Do the volumes of integer knotoids increase monotonically with the number of crossings? This must be true for a large enough crossing number, and the expectation is yes, but the question remains open. 
\end{question}

\section{Planar Knotoids} \label{planar_knotoid_section}

We now introduce our model for knotoid diagrams living in $\mathbb{R}^2$. As a reminder, knotoid diagrams in $\mathbb{R}^2$ differ from $S^2$ in the sense that strands of the knotoid cannot be pulled around the backside of the sphere in order to appear on the other side of the diagram. This restriction increases the number of distinct knotoids. Also, recall that we may consider planar knotoid diagrams either in $\mathbb{R}^2$ or $D^2$. 

\subsection{Maps from Planar Knotoids to Knots in Handlebodies}
We describe two maps, $\phi_{\mathbb{R}^2}^D$ and $\phi_{\mathbb{R}^2}^G$, taking planar knotoids to knots in a genus 3 handlebody $H_3$ and to knots in a genus 2 handlebody $H_2$, respectively. We will denote the set of knots in $H_3$ by $\mathcal{K}(H_3)$ and the set of knots in $H_2$ by $\mathcal{K}(H_2)$.  

Note that in \cite{BBHL}, the authors describe a map of planar knotoids to knots in the solid torus. First, they start with a planar knotoid diagram which they consider as a diagram in the disk. Then they take the double branched cover over the two endpoints to  obtain a knot projection on an annulus. Upon thickening, they obtain a knot in the solid torus. 
We give two alternate constructions that parallel our constructions in the case of knotoids in $S^2$.

\begin{figure}[htbp] 
\includegraphics[scale=0.4]{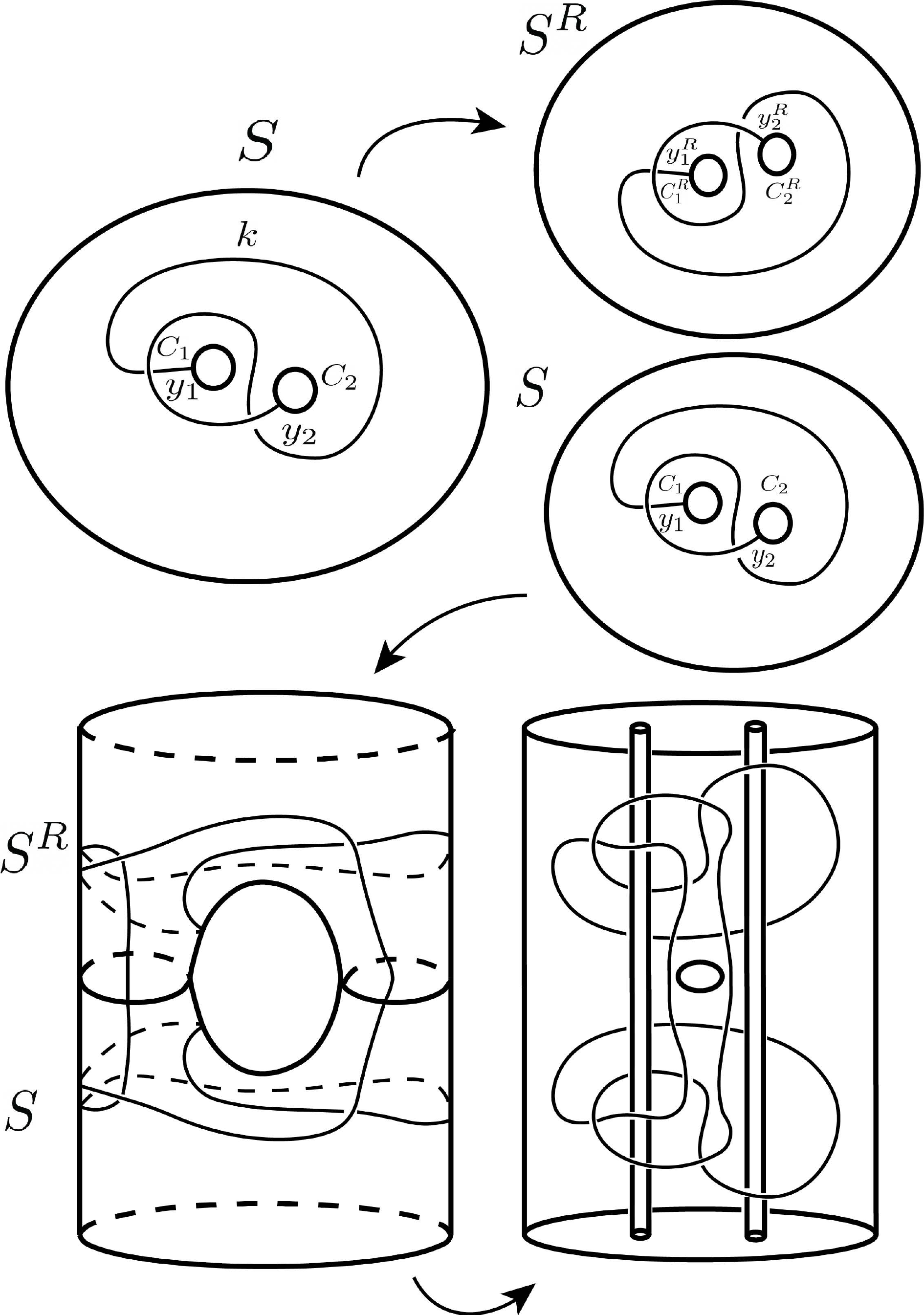}
\caption{The reflected doubling map $\phi_{\mathbb{R}^2}^D$ described in Definition \ref{defnPlanarDoublingMap}.}
\label{SphericalDoublingImage}
\end{figure}

\begin{definition} \label{defnPlanarDoublingMap}
Consider a knotoid diagram $\mathcal{D}$ on $D^2$ representing the planar knotoid $k$. Denote the two endpoints of $k$ by $x_1$ and $x_2$, and by abuse of notation, denote the immersed arc by $k$. Take disk neighborhoods $D_i$ of the $x_i$ for $i = 1, 2$, chosen sufficiently small so that the circle boundaries $C_i := \partial D_i$ are punctured exactly once by $k$. Denote the punctures by $y_i$. Consider the surface $S = D^2 \setminus (\mathring{D_1} \cup \mathring{D_2})$ and take a reflected copy $S^R$, which contains reflected copies $C_i^R$ once-punctured by $y_i^R$. Glue $S$ to $S^R$ via an orientation-preserving homeomorphism from $C_i$ to $C_i^R$ such that $y_i$ and $y_i^R$ are identified, for $i = 1, 2$. This yields a knot diagram on a twice-punctured torus, equivalently a knot $K$ in a thickened twice-punctured torus, which is homeomorphic to a genus 3 handlebody $H_3$. This defines a map $\phi_{\mathbb{R}^2}^D: \mathbb{K}(\mathbb{R}^2) \to \mathcal{K}(H_3)$ which we again call the \textit{reflected doubling map}. We write $\phi_{\mathbb{R}^2}^D(k) = K$ and refer to $K$ as the image of $k$ under the reflected doubling map. See Figure \ref{SphericalDoublingImage}. 
\end{definition}

\begin{figure}[htbp]
\includegraphics[scale=0.4]{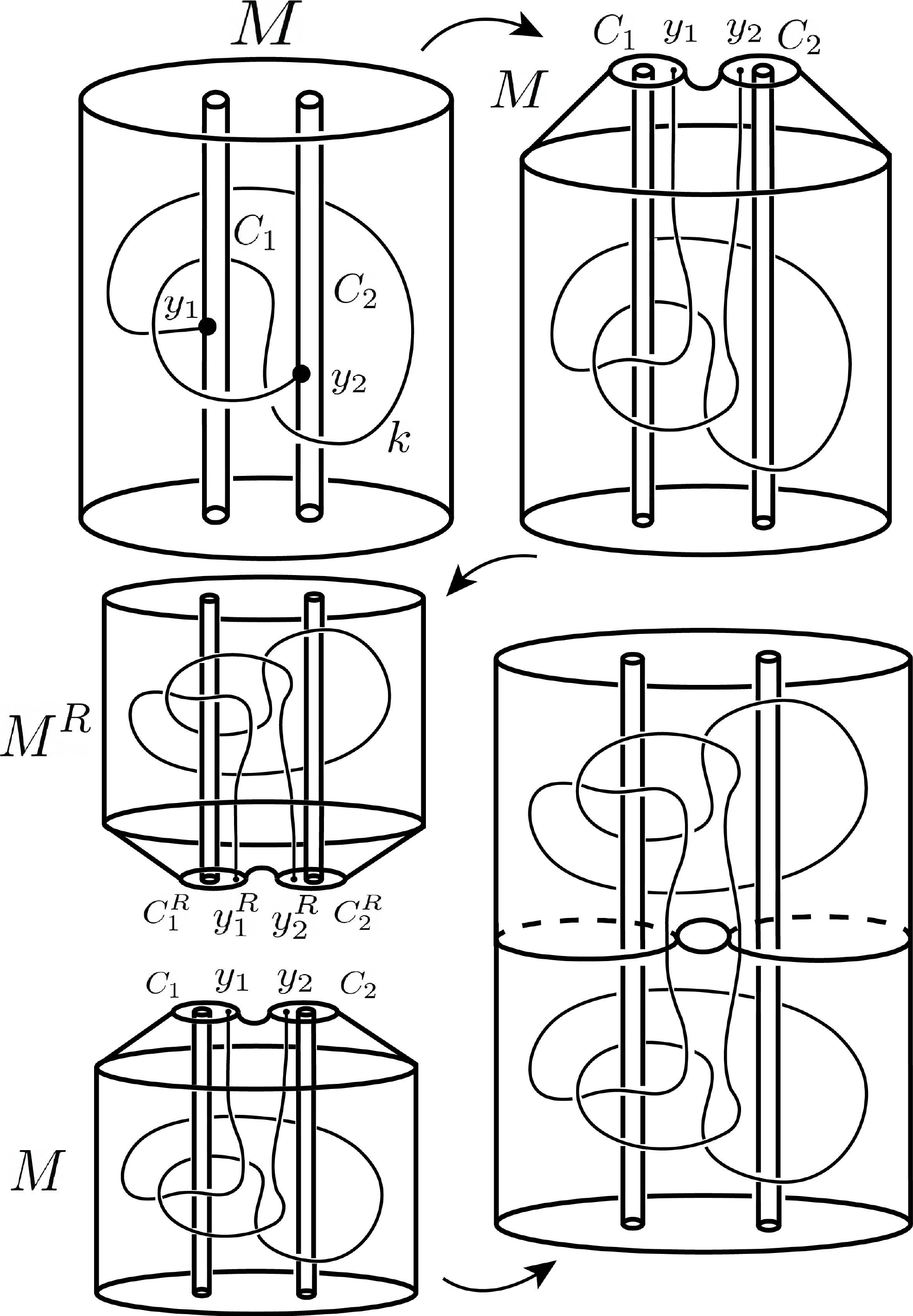}
\caption{An equivalent way to view the reflected doubling map $\phi^D_{\mathbb{R}^2}$ described in Definition \ref{defnPlanarDoublingMap}.}
\label{EquivalentPlanarDoublingFigure}
\end{figure}

Equivalently, we may visualize this construction in the planar rail diagram corresponding to $k$. That is, we start with a thickened disk $D^2 \times I$ together with two rails $\ell_i$ and an embedded arc $k$. Let $M$ be $(D^2 \times I) \setminus \mathring{N}(\ell_1 \cup \ell_2)$, which is just a thickened once-punctured torus. Let $C_i = \partial N(\ell_i)$ be the cylindrical boundaries of the rail neighborhoods, so that each is once-punctured by $k$. Take a reflected copy $M^R$, which contains reflected copies $C_i^R$ of the $C_i$. Then glue $M$ to $M^R$ along the $C_i$ via reflections, such that the punctures from $k$ are identified. We might visualize this as pulling the $C_i$ up to the top surface of $M$ before gluing in a ``natural" way. See Figure \ref{EquivalentPlanarDoublingFigure}. 

Note that we do not remove the boundary (compare to the spherical reflected doubling map) because the boundary is a genus 3 surface. In particular, we can put a hyperbolic metric on it. 

Now we may define a notion of hyperbolicity for planar knotoids based on this construction. 

\begin{definition} \label{defnPlanarDoublingHyperbolic}
Let $k \in \mathbb{K}(\mathbb{R}^2)$ be a planar knotoid, and let $K = \phi_{\mathbb{R}^2}^D(k)$ be the image of $k$ under the planar reflected doubling map. We say that $k$ is \textit{hyperbolic (under the reflected doubling map)} if $K$ is tg-hyperbolic in $H_3$. That is, $H_3 \setminus N(K)$ admits a complete hyperbolic metric such that $\partial H_3$, a genus 3 surface, is totally geodesic in this metric. If $k$ is hyperbolic, then we may define its \textit{hyperbolic volume under the reflected doubling map} $\vol_{\mathbb{R}^2}^D(k)$  as
$$
\vol_{\mathbb{R}^2}^D(k) := \frac{1}{2} \vol(H_3 \setminus N(K)). 
$$
\end{definition}

Note that we must require the stronger condition of tg-hyperbolicity on $Y \setminus N(\mathcal{G})$ in order to guarantee that it has a well-defined finite hyperbolic volume. Indeed, in general it is possible for a manifold with  boundary components of genus at least two to admit many hyperbolic metrics with varying (possibly infinite) volumes if we do not specify that the higher genus boundaries are totally geodesic in the metric. 

Similar to the spherical setting, we might choose to just glue $C_1$ to $C_2$ without doubling (compare to the spherical gluing map in Definition \ref{defnSphericalGluingMap}). This yields another map on planar knotoids, which we define now. 

\begin{definition} \label{defnPlanarGluingMap} 
Consider a knotoid diagram $\mathcal{D}$ on $D^2$ representing the planar knotoid $k$. We keep the setup and notation $x_i, D_i, C_i, S$, and $y_i$ as in Definition \ref{defnPlanarDoublingMap}. Instead of doubling, simply glue $C_1$ to $C_2$ via an orientation-reversing homeomorphism such that $y_1$ and $y_2$ are identified. This yields a knot diagram on a once-punctured torus which is equivalent to a knot $K$ in the genus 2 handlebody $H_2$. Then this defines a map $\phi_{\mathbb{R}^2}^G: \mathbb{K}(\mathbb{R}^2) \to \mathcal{K}(H_2)$ which we call the \textit{planar gluing map}. We write $\phi_{\mathbb{R}^2}^G(k) = K$ and refer to $K$ as the image of $k$ under the gluing map.

\begin{figure}[htbp]
\includegraphics[scale=0.4]{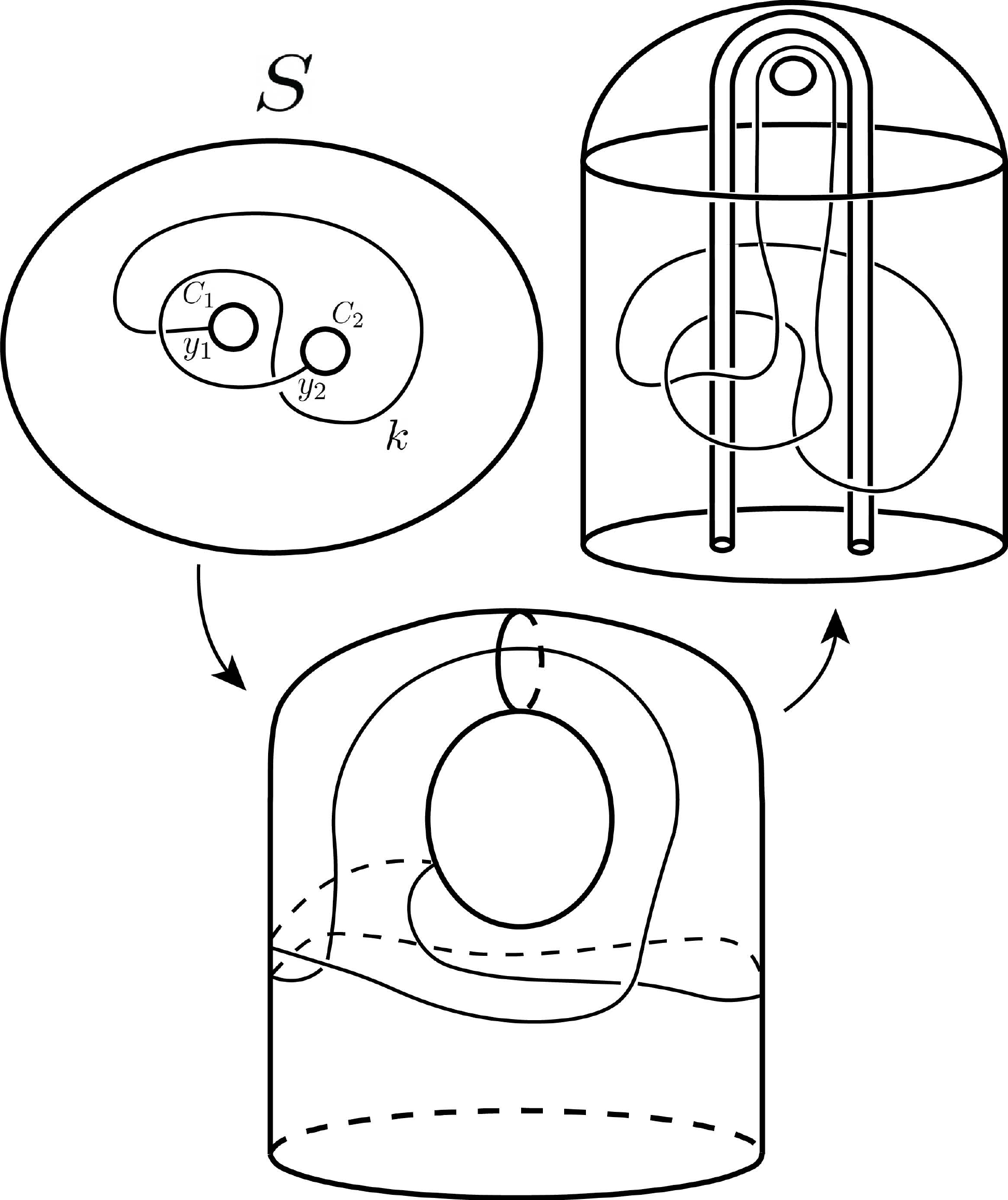}
\caption{The planar gluing map $\phi^G_{\mathbb{R}^2}$ described in Definition  \ref{defnPlanarGluingMap}.}
\label{PlanarGluingImage}
\end{figure}

\end{definition}

Equivalently, we may visualize this construction in the planar rail diagram corresponding to $k$. The steps are the same as for the planar reflected doubling map, except we now glue the cylindrical boundaries $\partial N(\ell_i)$ of the rail neighborhoods to each other in the ``natural" way. The gluing map gives us another definition of hyperbolicity for planar knotoids, which we define now.

\begin{figure}
\includegraphics[scale=0.4]{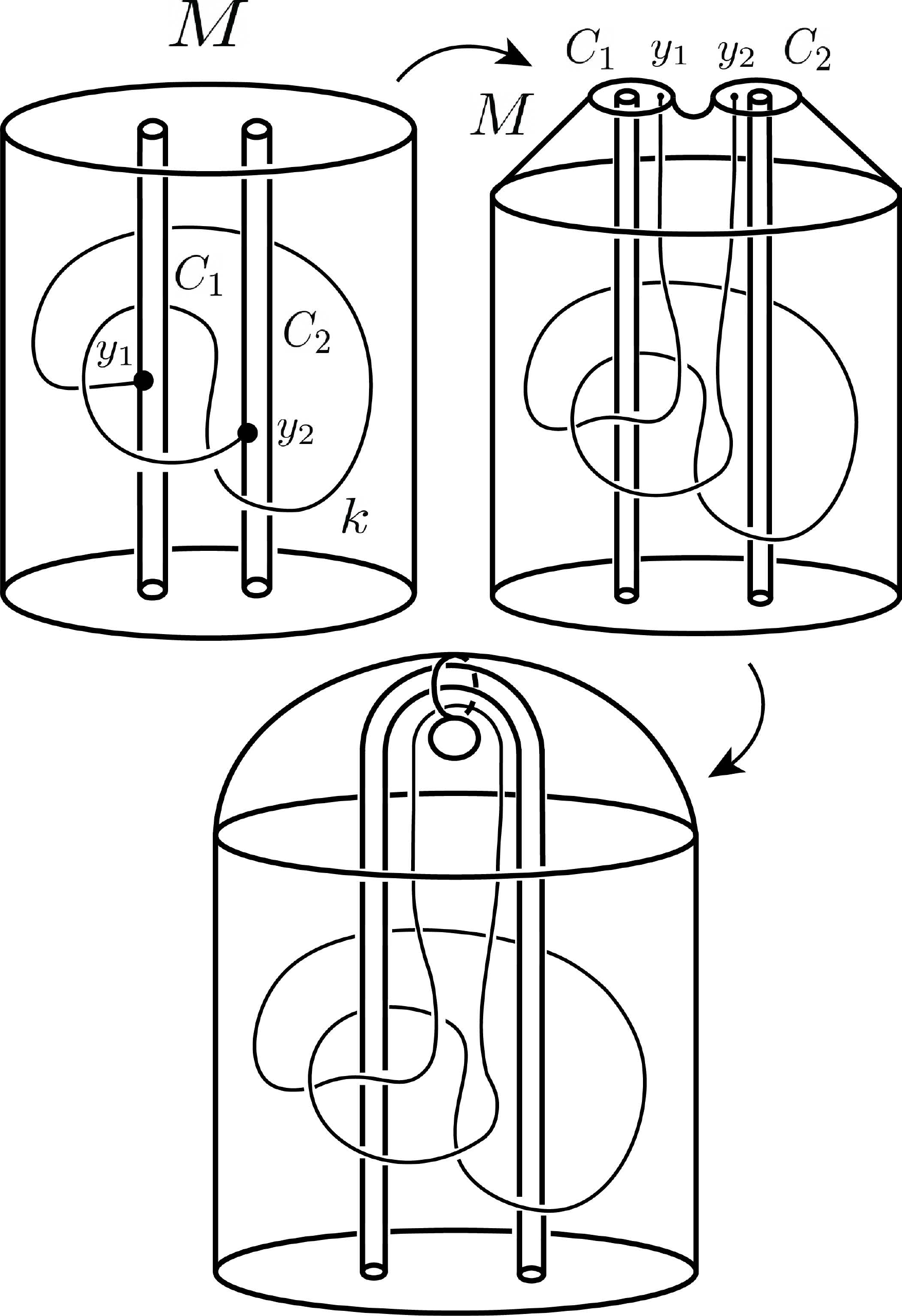}
\label{EquivalentPlanarGluing}
\caption{An equivalent way to view the planar gluing map $\phi^G_{\mathbb{R}^2}$ from Definition \ref{defnPlanarGluingMap}.}
\end{figure}

\begin{definition}
Let $k \in \mathbb{K}(\mathbb{R}^2)$ be a planar knotoid, and let $K = \phi_{\mathbb{R}^2}^G(k)$ be its image under the planar gluing map. We say that $k$ is \textit{hyperbolic (under the gluing map)} if $K$ is tg-hyperbolic in $H_2$. That is, $H_2 \setminus N(K)$ admits a complete hyperbolic metric such that $\partial H_2$, a genus 2 surface, is totally geodesic in this metric. If $k$ is hyperbolic, then we may define its \textit{hyperbolic volume under the gluing map} $\vol_{\mathbb{R}^2}^G(k)$ as
$$
\vol_{\mathbb{R}^2}^G(k) := \vol(H_2 \setminus N(K)). 
$$
\end{definition}

It is natural to ask how hyperbolicity for planar knotoids under the reflected doubling map and the gluing map are related, and how the volumes compare. For the first question, in contrast with the spherical setting, being hyperbolic under the reflected doubling map implies hyperbolic under the gluing map, but not the converse. For the second question, it can be shown that the volume under the gluing map cannot be less than the volume under the reflected doubling map. 

\begin{proposition}
For a knotoid $k$, if $\phi_{\mathbb{R}^2}^D(k)$ is tg-hyperbolic, then $\phi_{\mathbb{R}^2}^G(k)$ is tg-hyperbolic. That is, tg-hyperbolicity in the reflected doubling map implies tg-hyperbolicity in the gluing map. Furthermore, $\operatorname{vol}(\phi_{\mathbb{R}^2}^G(k)) \geq \operatorname{vol}(\phi_{\mathbb{R}^2}^D(k))$. 
\end{proposition}

\begin{proof}

Let $M$ be the manifold $H_2\setminus \mathring{N}(\phi^G_{\mathbb{R}^2})$, and let $F \subset M$ be the shaded surface as in Figure \ref{FinM}. 

\begin{figure}
    \centering
    \includegraphics[scale=0.45]{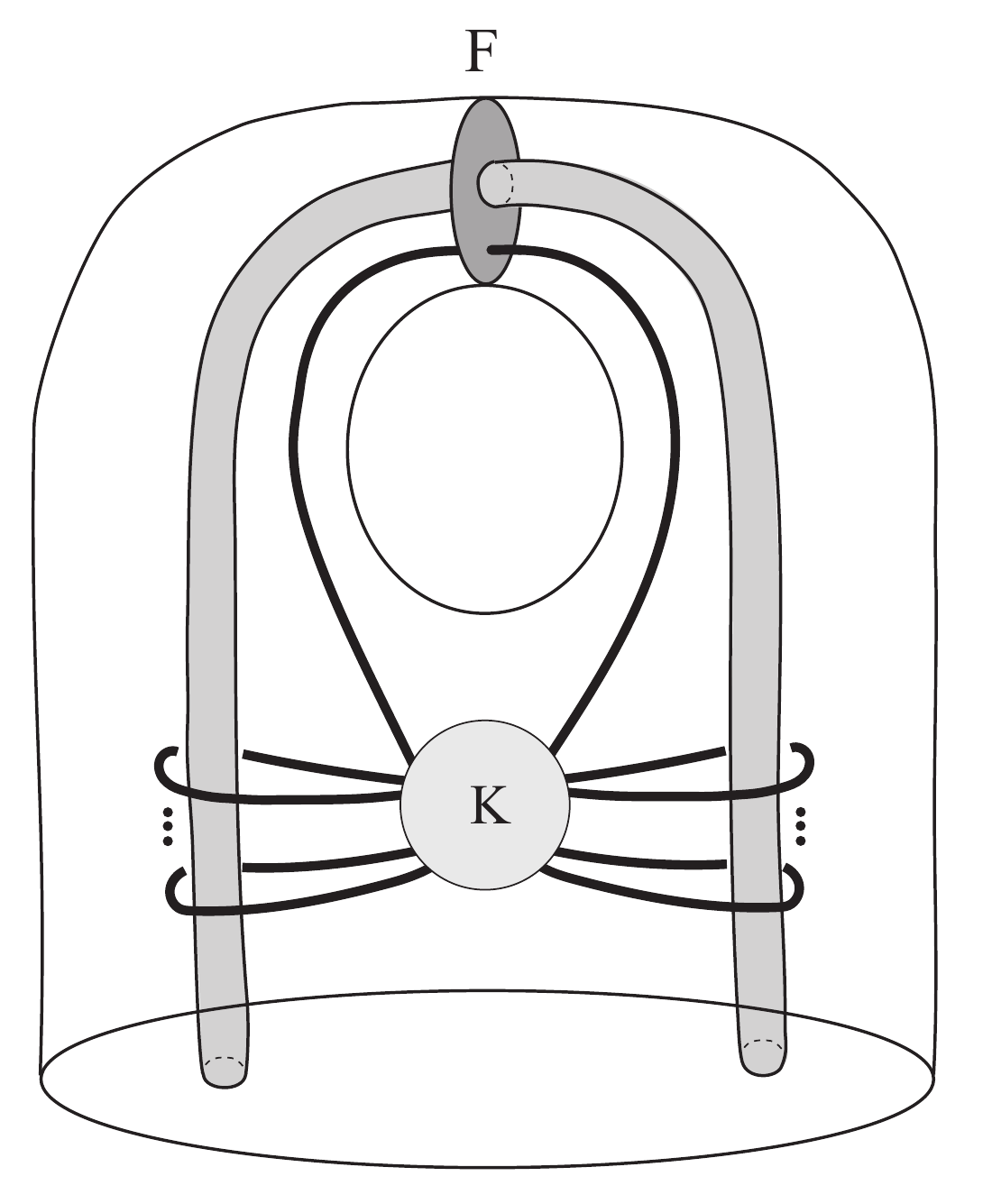}
    \caption{The punctured annulus $F$ sits in $M$.}
    \label{FinM}
\end{figure}

To obtain the double cover of $M$, we cut $M$ along $F$ and take a second rotated copy that we  glue to the original copy along the copies of $F$, as shown in Figure \ref{MandM'}(a). Note that $M$ is hyperbolic if and only if $M'$ is also. 

To obtain the reflected double, we again cut $M$ open along $F$, but now we take a reflected copy of the result, which we glue along the copies of $F$ to obtain the manifold $M''$ as in Figure \ref{MandM'}(b). Note that $F''$ is totally geodesic as it is the fixed point set of the reflection of the manifold that results.

\begin{figure}
    \centering
    \includegraphics[scale=0.4]{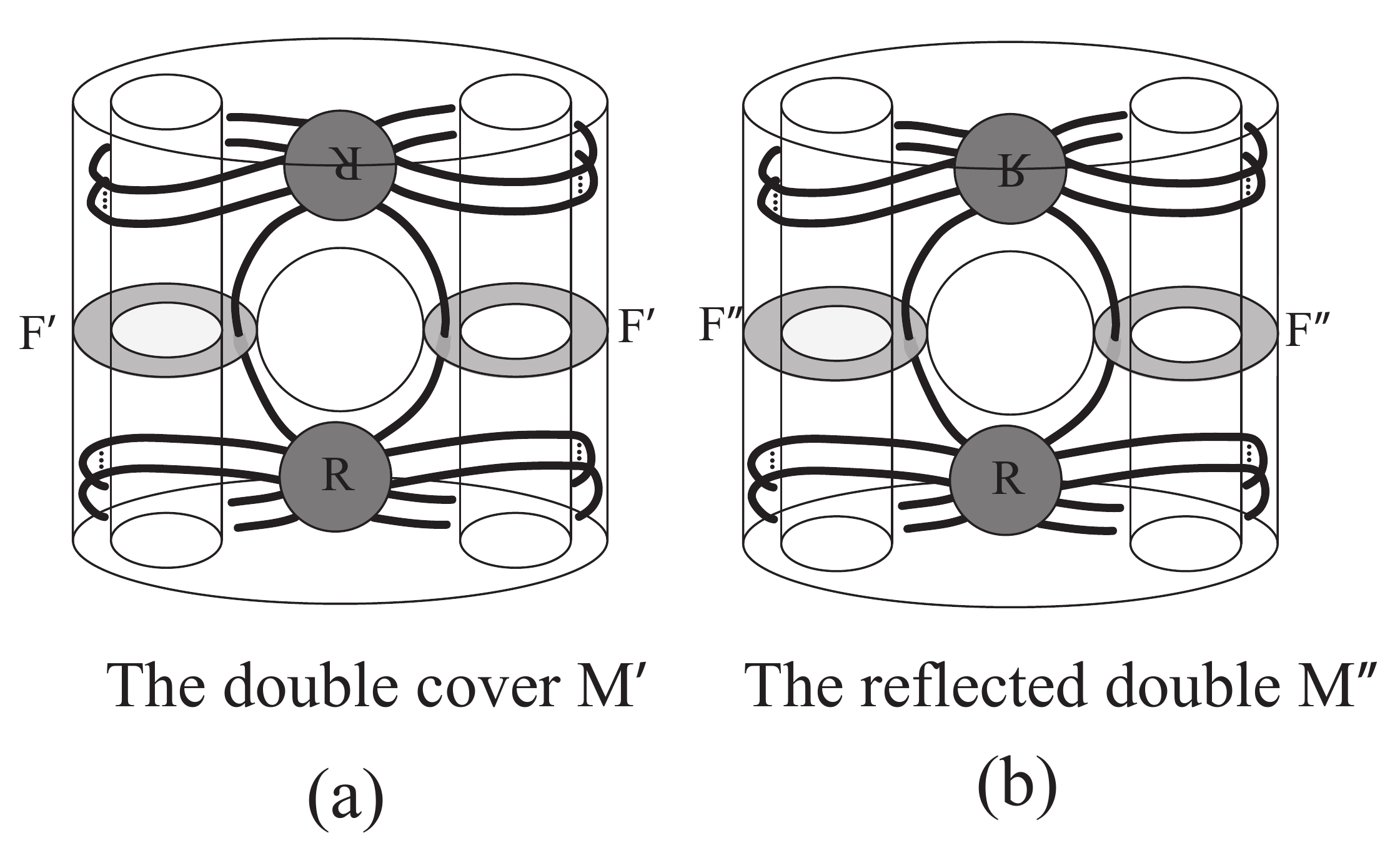}
    \caption{ To obtain $M'$ from $M''$, cut along $F'' \subset M''$, reflect the top component and reglue along the copies of $F''$ as shown in the diagram.}
    \label{MandM'}
\end{figure}

To obtain the double cover $M'$  of $M$ from the reflected double $M''$, we cut along the totally geodesic surface $F''$, reflect the top component along a vertical plane perpendicular to the page and glue back together along the copies of $F''$. Results from  \cite{agol_storm_thurston} and its extension in \cite{calegari_freedman_walker} tell us that the resulting manifold $M'$ is also hyperbolic, and that the volume of $M'$ is at least as large as the  volume of the reflected double  $M''$. Furthermore, the double cover $M'$ is hyperbolic if and only if $M$ is hyperbolic, and its volume is twice the volume of the double cover. Since by definition, we halve the volume of $M'$ to obtain the volume of a knotoid as in \ref{defnPlanarDoublingMap}, the claim follows. 

\end{proof}

\begin{figure}
    \centering
    \includegraphics[scale = 1.2]{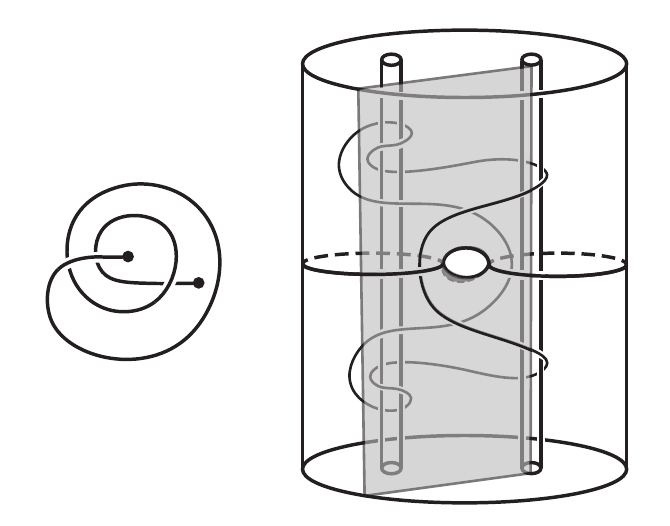}
    \caption{The planar knotoid $3_4$ and its image under $\phi_{\mathbb{R}^2}^D$ in $H_3$ with an essential annulus shaded.}
    \label{planarknotoid34}
\end{figure}

However, the converse is not true. For example, both of the planar knotoids $2_2$ and $3_4$ (using the notation of \cite{tabulation}) are tg-hyperbolic under the gluing construction, with volumes of $9.13447\dots$ and $10.68436\dots$, respectively. However, neither knotoid is hyperbolic under the doubling construction. In Figure \ref{planarknotoid34}, when $k$ is the planar knotoid $3_4$, we exhibit an essential annulus in the manifold $H_3 \setminus \mathring{N}(\phi^D_{\mathbb{R}^2}(k))$  that obstructs hyperbolicity. 

\subsection{An Infinite Family of Hyperbolic Planar Knotoids}
 
\begin{figure}[htbp]
\centering
\includegraphics[scale=0.6]{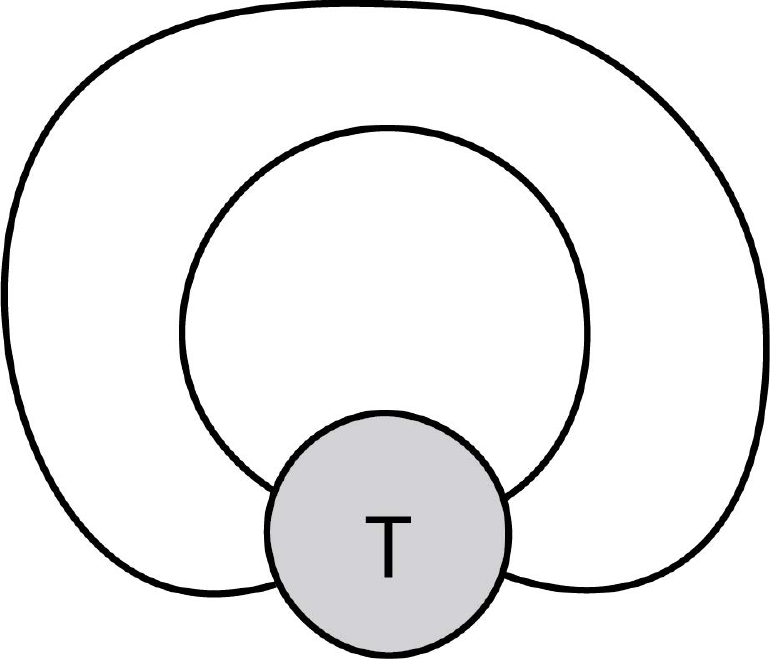}
\hspace{1cm}
\includegraphics[scale=0.6]{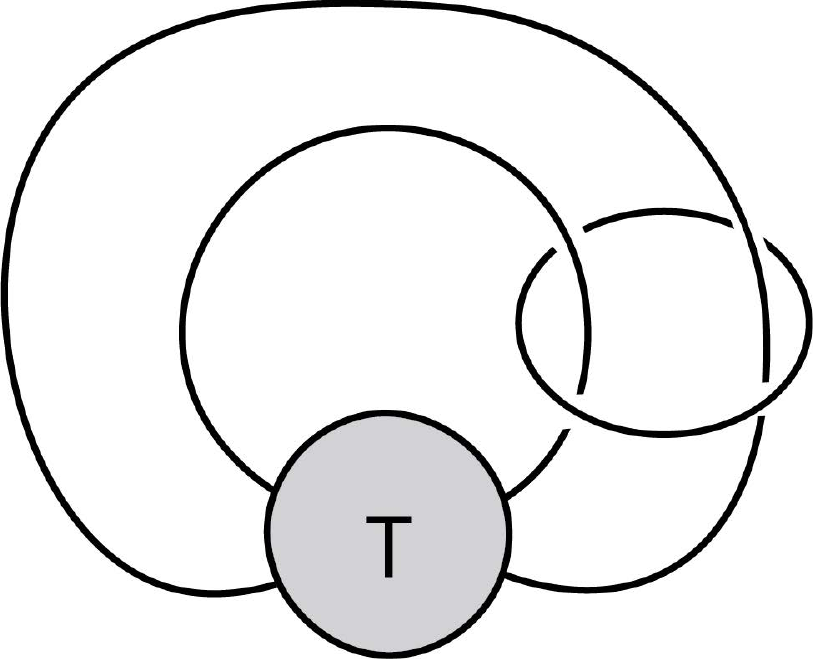}
\caption{The link $L_T$ and its augmentation $\overline{L}_T$}
\label{LT_and_augmentedLT}
\end{figure}

Using the results of \cite{CompositionPaper}, we can construct an infinite family of knotoids that are hyperbolic under the map $\phi^G_{\mathbb{R}^2}$. For a tangle $T$, we denote by $L_T$ the link which is the closure of $T$ as shown in Figure \ref{LT_and_augmentedLT}, and we denote the augmentation of $L_T$ by $\overline{L}_T$.

\begin{figure}[htbp]
 \centering
 \includegraphics[scale=0.55]{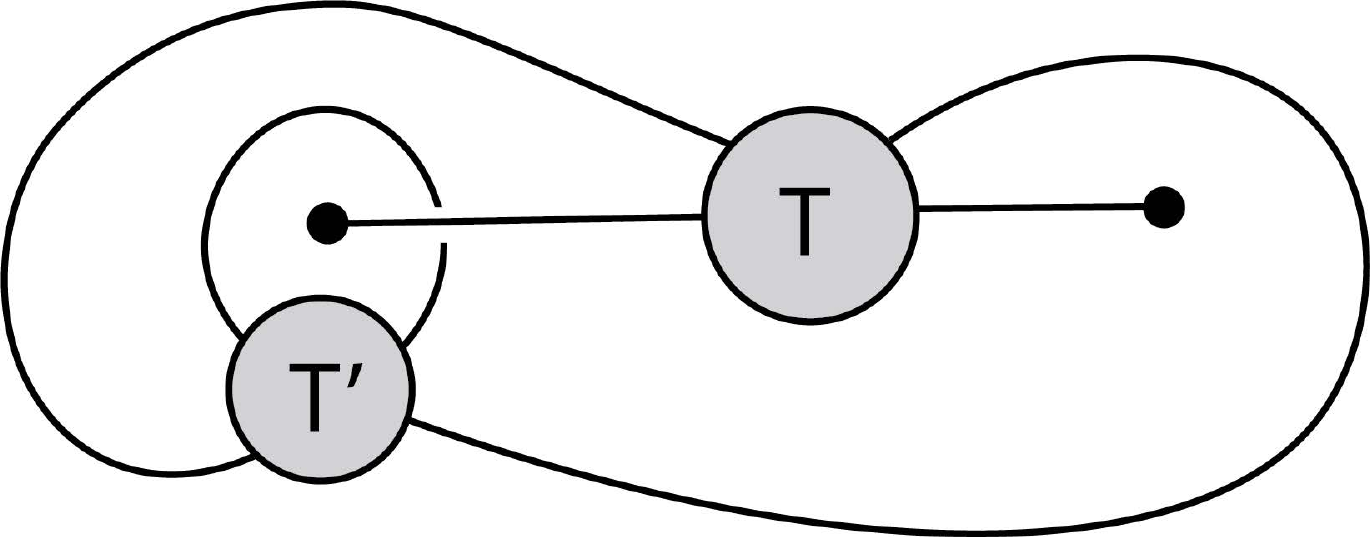}
 \caption{The planar knotoid $k_{T,T'}$ associated to $T,T'$}
 \label{Infinite_family_knotoid}
 \end{figure}

\begin{theorem}Let $T,T'$ be two tangles such that the links $\overline{L}_T,\overline{L}_{T'}$  are hyperbolic links with two components. Then $\phi_{\mathbb{R}^2}^G(k_{T,T'})$ is $tg$-hyperbolic, where $k_{T,T'}$ is the knotoid shown in Figure \ref{Infinite_family_knotoid}. 
\label{Infinite_family_Theorem}
\end{theorem}

In particular, by the results of \cite{colin!} the hypotheses hold if the tangles $T,T'$ are prime and such that $L_T,L_{T'}$ are alternating and not $2$-braids.

\section{Volume Computations} \label{volume_tables}

All volume computations in this section were carried out using the software SnapPy \cite{SnapPy}.
\subsection{Spherical Knotoids}

Of all of the knotoids in $S^2$ with 5 or fewer crossings tabulated in \cite{tabulation}, of which there are 36, only five are not hyperbolic. Moreover, all five of these non-hyperbolic knotoids are knot-type, a class we already know to be non-hyperbolic as in Proposition \ref{nonhyperbolicity_of_knottype}.  This data suggests that similarly to classical knots and to virtual knots, hyperbolic knotoids dominate non-hyperbolic knotoids for small crossing number. 

In the tables below, we include diagrams for the first several knotoids; the notation for the other knotoids comes from the tables of knotoids in \cite{tabulation}. 

\medskip

\begin{center}
    \begin{tabular}{| c | c | c | c |} \hline
 Knotoid & Volume & Knotoid & Volume \\  
 \hline
 \includegraphics[scale=0.2]{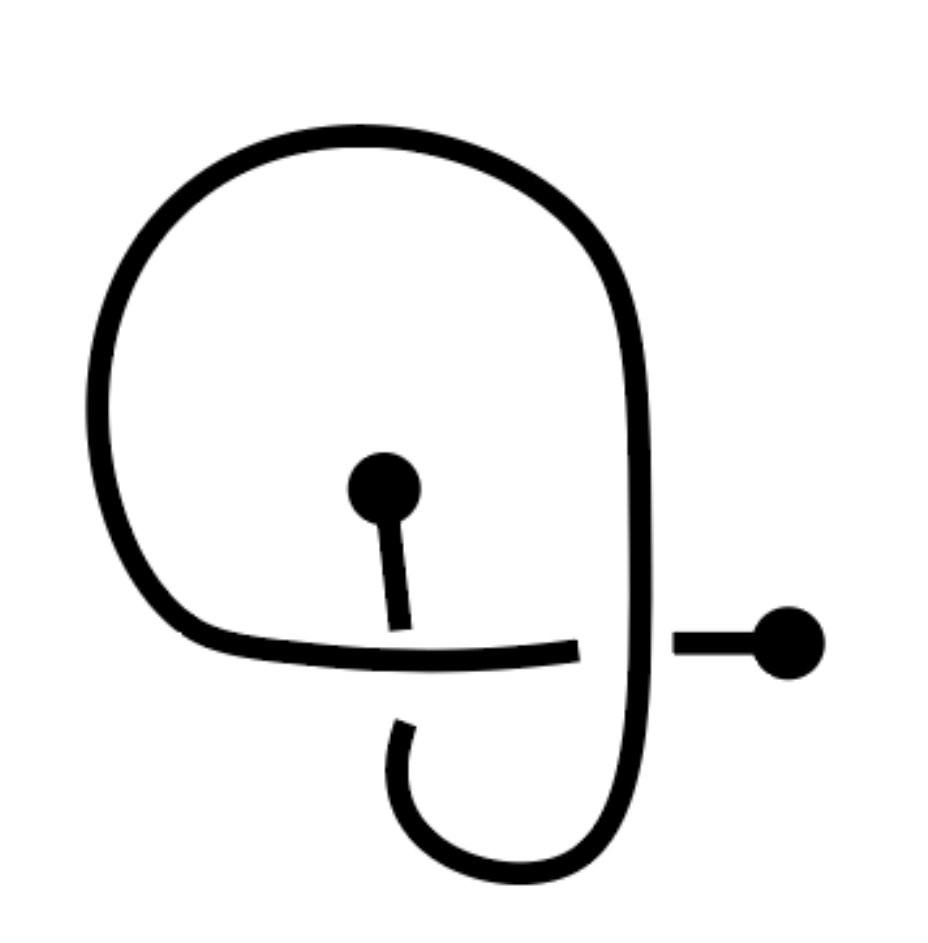} $2_1$ & $5.33349\dots$ &  \includegraphics[scale=0.2]{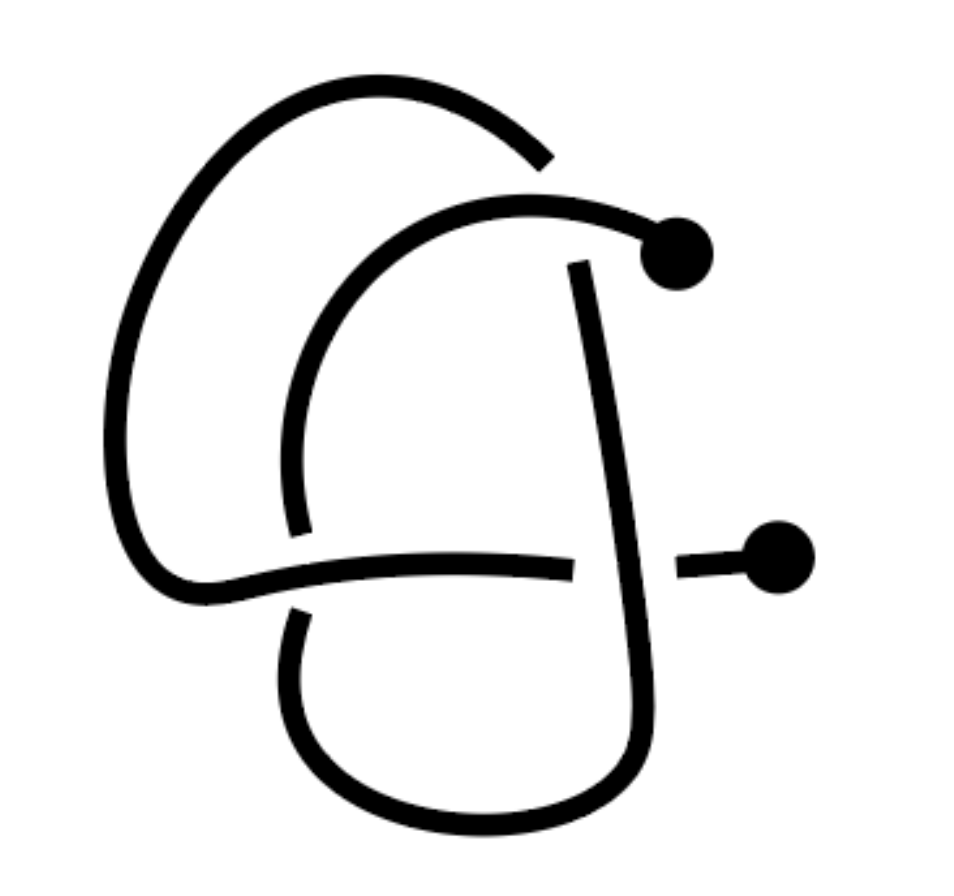} $3_1$ & Non-hyperbolic\\ 
 \hline
 
\includegraphics[scale=0.2]{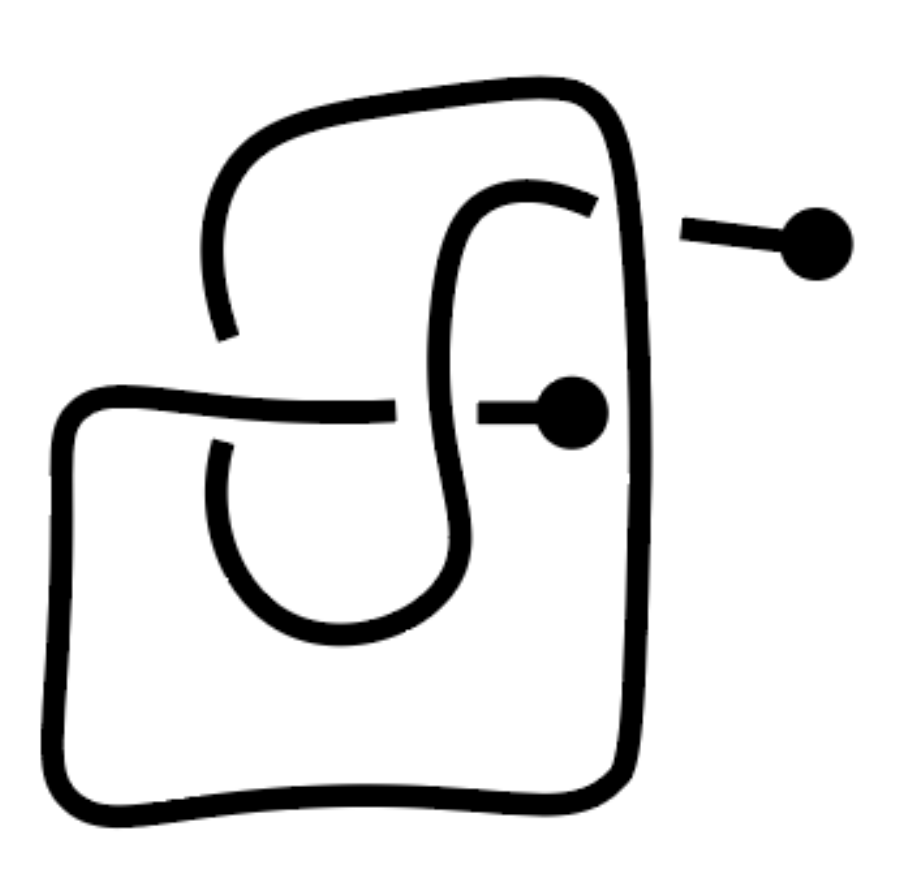} $3_2$ & $7.70691\dots$ & \includegraphics[scale=0.2]{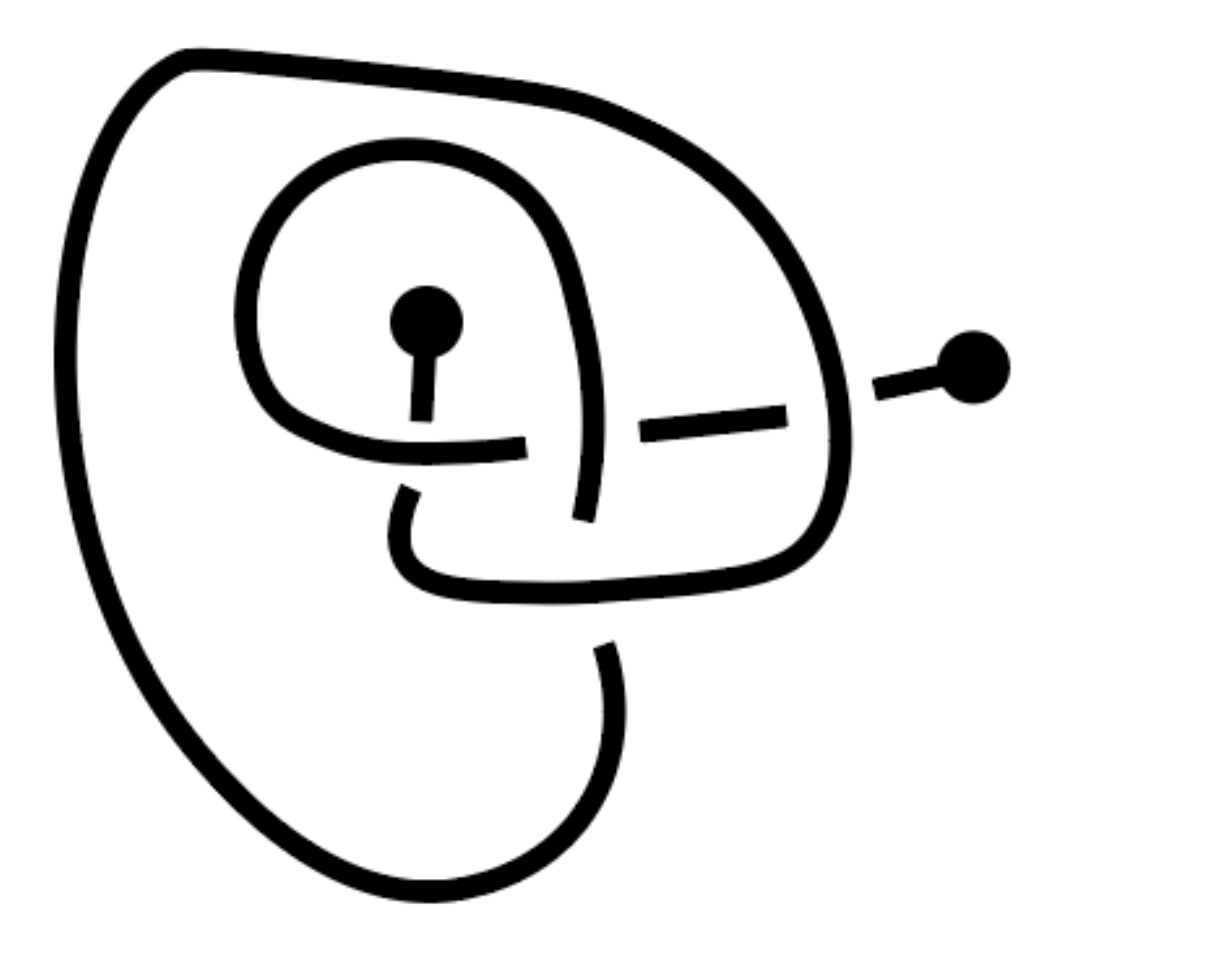} $4_2$ & $8.79335\dots$\\ 
\hline

\includegraphics[scale=0.2]{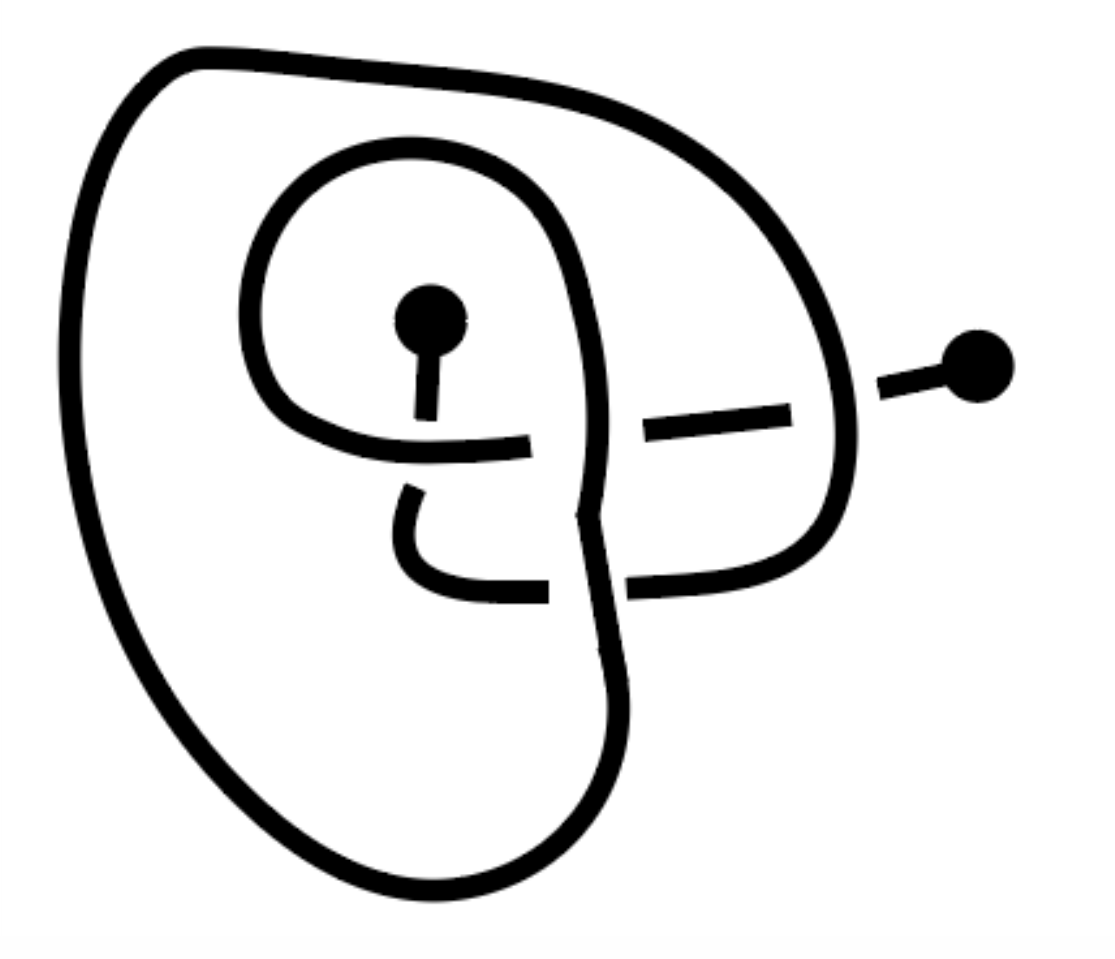} $4_3$ & $6.55174\dots$ & \includegraphics[scale=0.2]{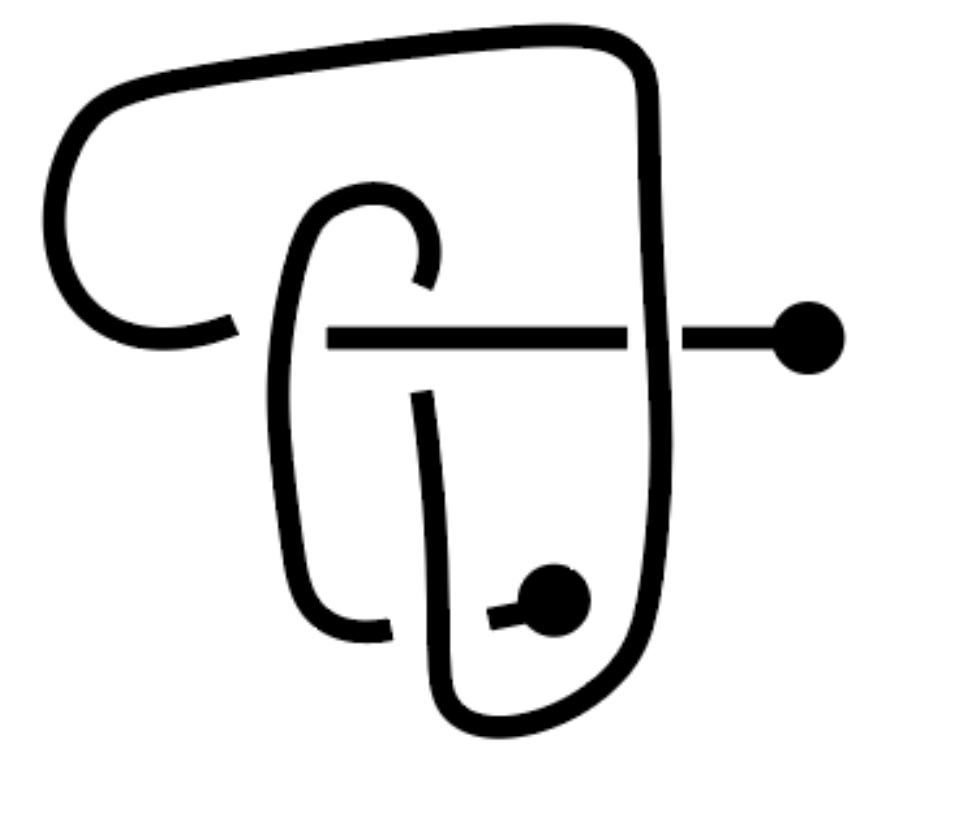} $4_4$ & $9.96651\dots$\\
\hline

\includegraphics[scale=0.2]{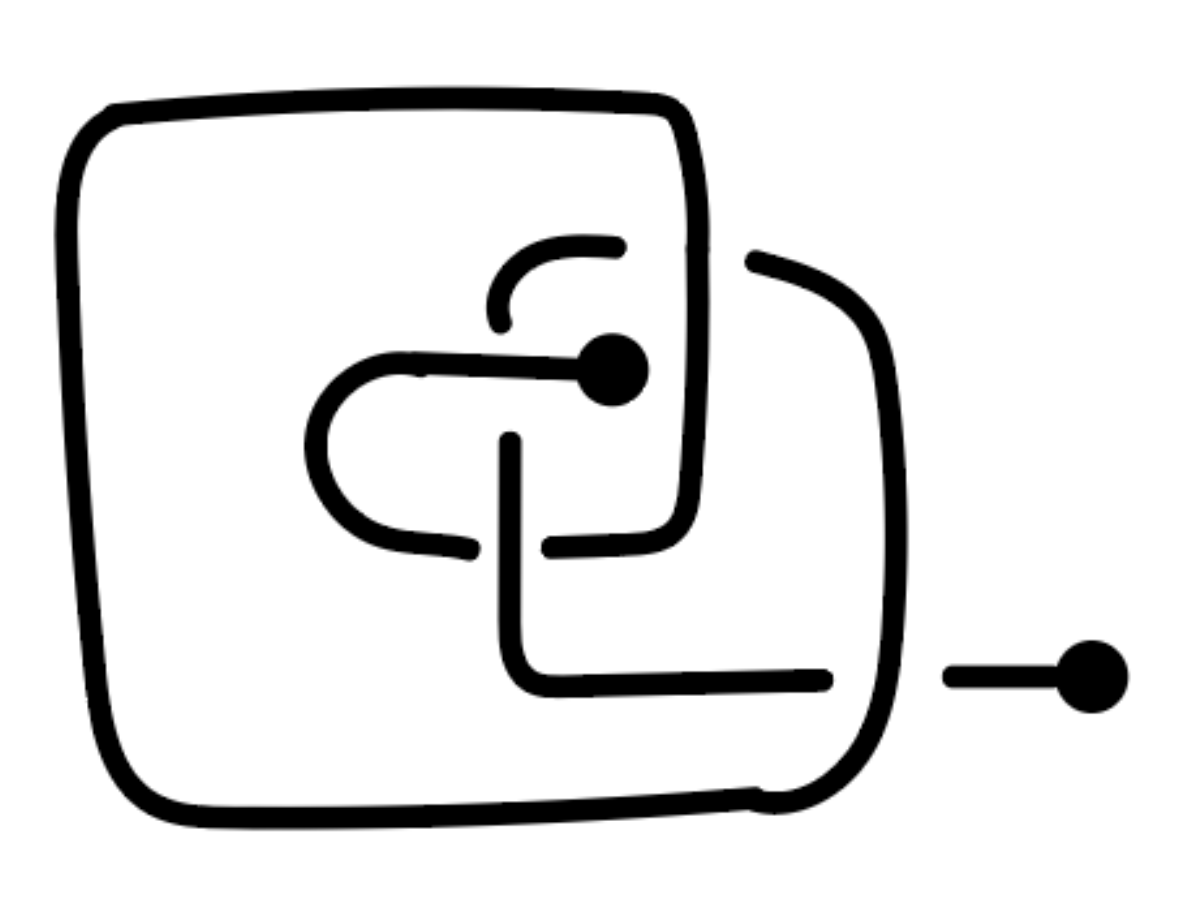} $4_5$ & $10.55687\dots$ & \includegraphics[scale=0.2]{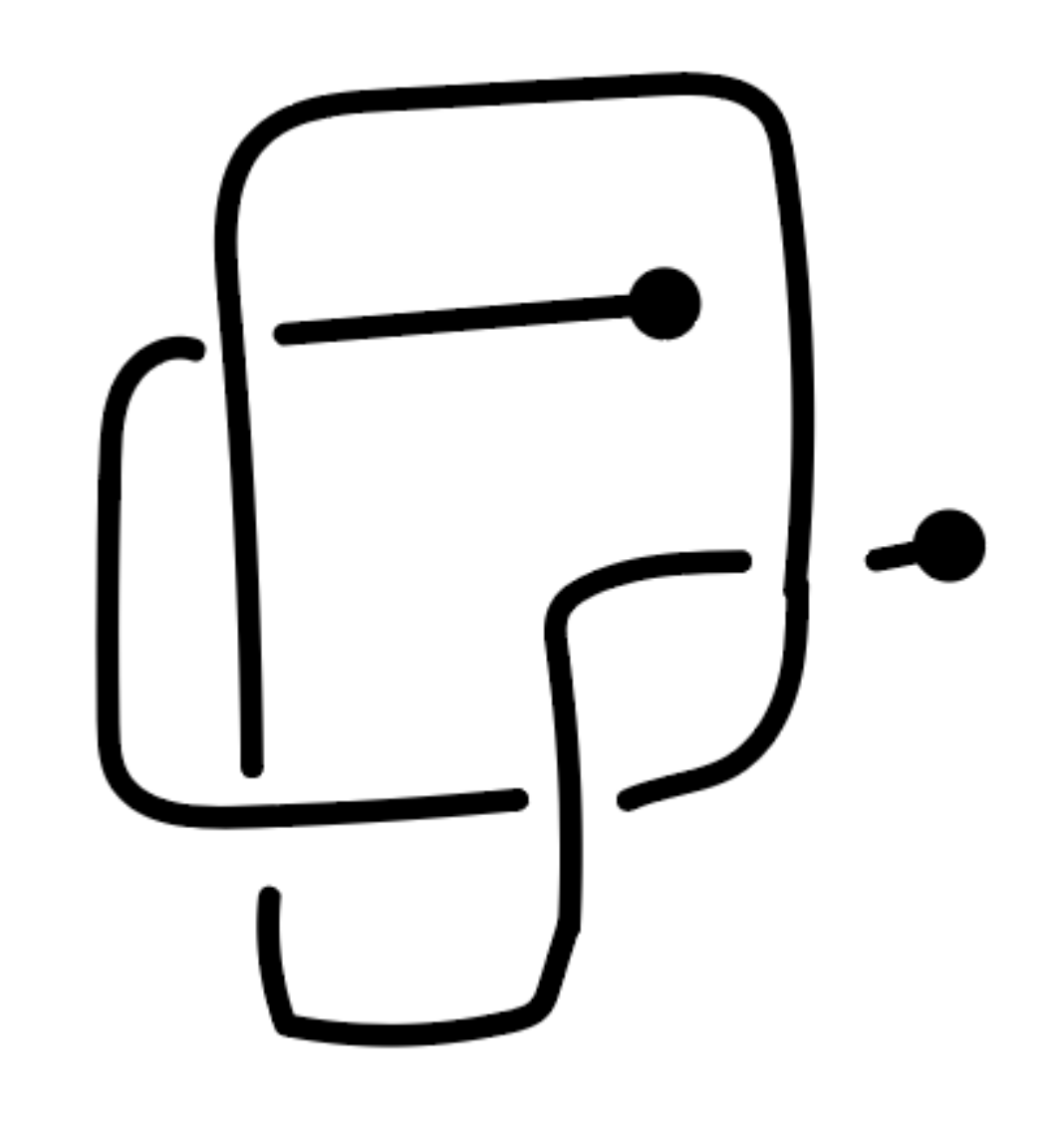} $4_6$ & $8.35550\dots$\\
\hline

\includegraphics[scale=0.2]{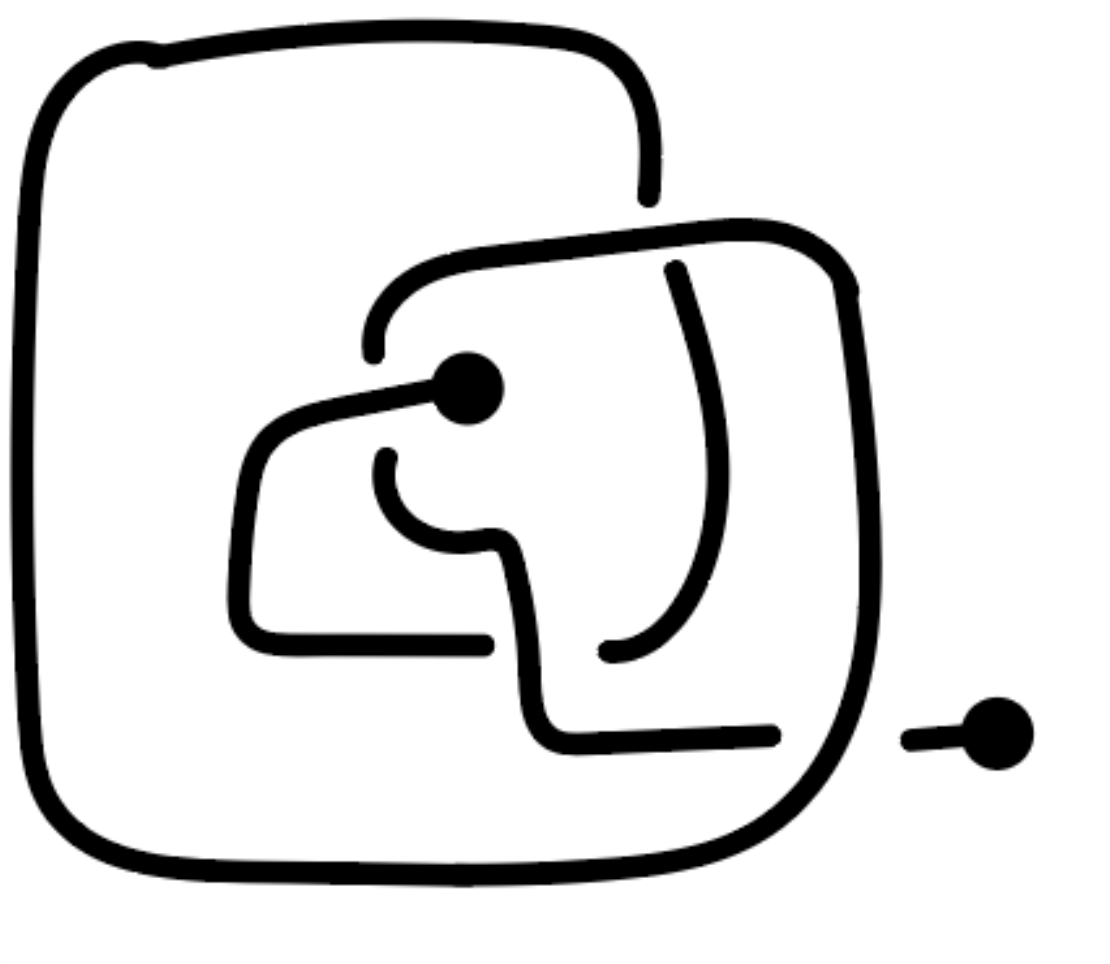} $4_7$ & $8.92932\dots$ & \includegraphics[scale=0.2]{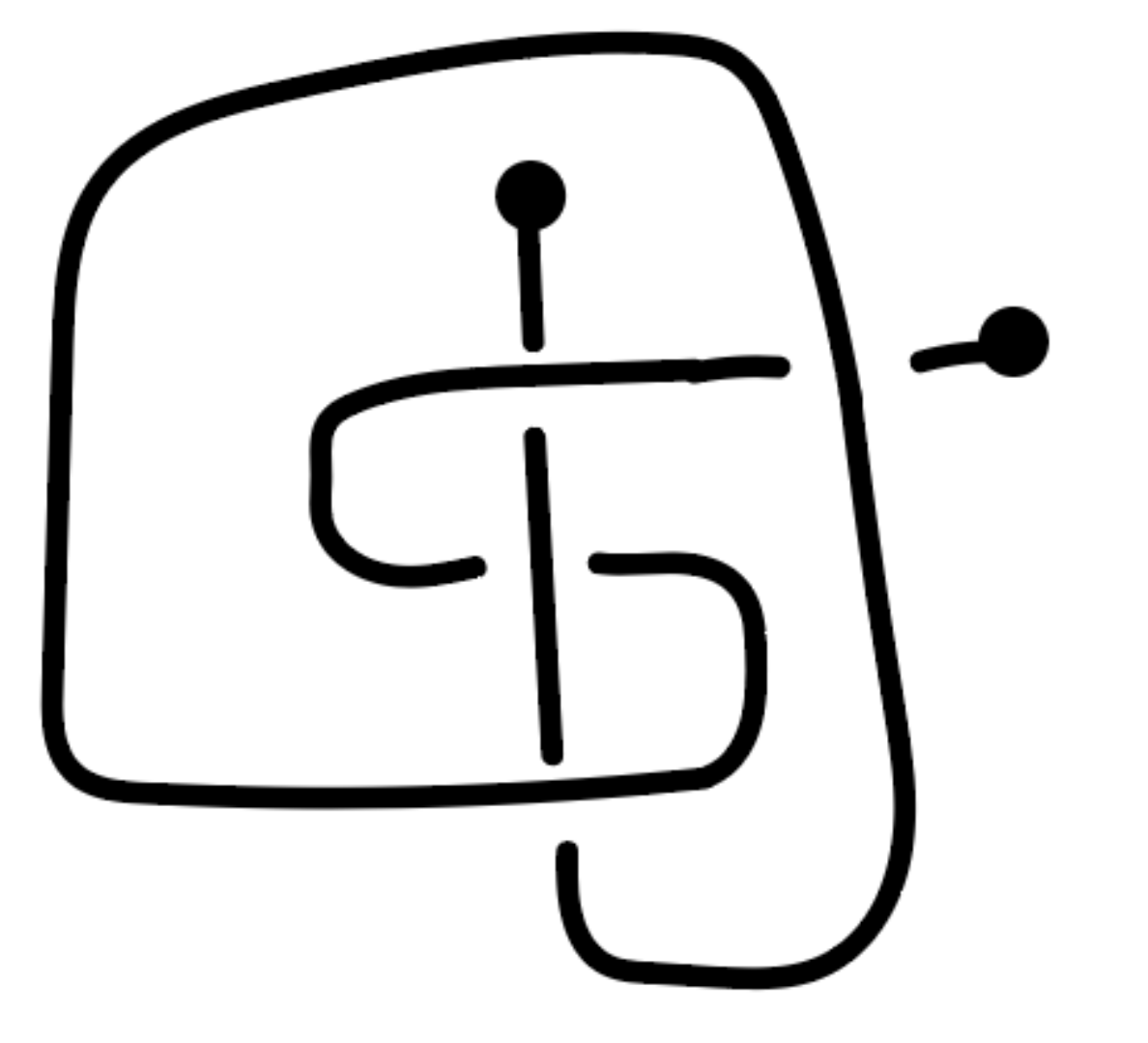} $4_8$ & $8.96736\dots$\\
\hline

$5_3$ & $11.76223\dots$ & $5_4$ & $10.74026\dots$ \\
\hline

$5_5$ & $9.31234\dots$ & $5_6$ & $10.56281\dots$ \\
\hline

$5_7$ & $11.08217\dots$ & $5_8$ & $12.01109\dots$ \\ 
\hline

$5_9$ & $13.31846\dots$ & $5_{10}$ & $12.00595\dots$\\
\hline

$5_{11}$ & $12.27656\dots$ & $5_{12}$ & $12.21631\dots$\\
\hline

$5_{13}$ & $11.28460\dots$ & $5_{14}$ & $12.73265\dots$\\
\hline

$5_{15}$ & $14.61681\dots$ & $5_{16}$ & $13.4431\dots$ \\ 
\hline

$5_{17}$ & $13.54793\dots$ & $5_{18}$ & $10.39689\dots$ \\
\hline

$5_{19}$ & $12.64232\dots$ & $5_{20}$ & $13.99859\dots$ \\
\hline

$5_{21}$ & $13.70176\dots$ & $5_{22}$ & $12.54144\dots$\\
\hline 

$5_{23}$ & $13.43713\dots$ & $5_{24}$ & $16.76577\dots$\\
\hline

\end{tabular}
\end{center}

\subsection{Planar Knotoids}

 Interestingly, fewer planar knotoids are hyperbolic compared to spherical knotoids. As in the spherical case, any knotoid $k$ with a diagram $\mathcal{D}$ that is knot-type cannot be hyperbolic, but furthermore, if either endpoint of $k$ is in the exterior region of $\mathcal{D}$, then $k$ cannot be hyperbolic. In the table below, knotoids listed before $4_1$ in \cite{tabulation} and not included in this table are not hyperbolic in either construction.  

\begin{center}
    \begin{tabular}{|c | c | c |} \hline
    
    Knotoid & Volume in the Gluing Map & Volume in the Reflected Doubling Map \\
    \hline
    
    $2_2$ & $9.13447\dots$ & Non-hyperbolic\\
    \hline
    $3_3$ & $13.10319\dots$ & $10.57083\dots$\\
    \hline
    $3_4$ & $10.68436\dots$ & Non-hyperbolic \\
    \hline
    $3_5$ & $10.68436\dots$ & $10.15992\dots$ \\
    \hline
    $3_6$ & $13.10319\dots$ & $12.04805\dots$ \\
    \hline
    $3_8$ & $14.92611\dots$ & Non-hyperbolic \\
    \hline
    $3_{10}$ & $11.15142\dots$ & Non-hyperbolic \\
    \hline
    $3_{13}$ & $14.32335\dots$ & $14.21299\dots$ \\
    \hline
    $3_{14}$ & $14.10221\dots$ & $13.80027\dots$ \\
    \hline
    $3_{17}$ & $11.07797\dots$ & $10.36434\dots$ \\
    \hline
    $3_{19}$ & $12.37889\dots$ & $12.11410\dots$ \\
    \hline
    $3_{20}$ & $12.37889\dots$ & $11.95559\dots$ \\
    \hline

    \end{tabular}
\end{center}


\bibliographystyle{plain}
\bibliography{refs}
\nocite{*}
\end{document}